\documentclass[a4paper,11pt]{amsart}
\usepackage{amssymb,amsthm,bbm,epic,eepic,graphics,amsbsy}
\usepackage{longtable}
\usepackage{a4wide}
\usepackage{amssymb,amsthm}
\usepackage{xypic}
\usepackage{array}
\usepackage[applemac]{inputenc}
\newcommand{\la}{\lambda}

\DeclareMathOperator{\Irr}{\mathrm{Irr}}

\def\Ker{\mathrm{Ker}\,}
\def\Res{\mathrm{Res}}
\def\Ind{\mathrm{Ind}}
\def\Aut{\mathrm{Aut}}

\def\C{\ensuremath{\mathbbm{C}}}

\def\Z{\mathbbm{Z}}
\def\N{\mathbbm{N}}

\def\bar{\overline}

\def\eps{\epsilon}

\def\onto{\twoheadrightarrow}
\def\into{\hookrightarrow}

\def\GL{\mathrm{GL}}

\def\lcm{\mathrm{lcm}}
\def\Aff{\mathrm{Aff}}
\def\Id{\mathrm{Id}}
\def\OO{\mathcal{O}}
\newtheorem{theo}{Theorem}[section]
\newtheorem{theointr}{Theorem}
\newtheorem{fact}{Fact}
\newtheorem{prop}[theo]{Proposition}

\newtheorem{lemma}[theo]{Lemma}

\newtheorem{cor}[theo]{Corollary}

\def\mod{\ \mathrm{mod}\ }

\makeatletter
\DeclareRobustCommand\widecheck[1]{{\mathpalette\@widecheck{#1}}}
\def\@widecheck#1#2{%
   \setbox\z@\hbox{\m@th$#1#2$}%
   \setbox\tw@\hbox{\m@th$#1%
      \widehat{%
         \vrule\@width\z@\@height\ht\z@
         \vrule\@height\z@\@width\wd\z@}$}%
   \dp\tw@-\ht\z@
   \@tempdima\ht\z@ \advance\@tempdima2\ht\tw@ \divide\@tempdima\thr@@
   \setbox\tw@\hbox{%
      \raise\@tempdima\hbox{\scalebox{1}[-1]{\lower\@tempdima\box\tw@}}}%
   {\ooalign{\box\tw@ \cr \box\z@}}}
\makeatother

\def\WT{\widetilde{W}}
\def\WTH{\widetilde{W}_r}
\def\WH{W_r}
\def\Ad{\mathrm{Ad}}

\title{{\bf Branching properties for the groups $G\lowercase{(de,e,r)}$}}
\author{Ivan Marin}
\date{August 16th, 2007}
\begin{document}

\maketitle

\bigskip

\bigskip

\noindent {\bf Abstract.} We study general properties of
the restriction of the representations of
the finite complex reflection groupes $G(de,e,r+1)$
to their maximal parabolic subgroups of type $G(de,e,r)$,
and focus notably on the multiplicity of components.
In combinatorial terms, this amounts to the following
question : which symmetries arise or disappear
when one changes (exactly) one pearl in a combinatorial
necklace ?

\medskip


\section{Introduction}

\subsection{Motivations}

It is well-known that, for irreducible, classical Coxeter groups of
type $A_{n+1}$, $B_{n+1}$ and $D_{n+1}$, the restriction
of irreducible representations to their natural maximal parabolic
subgroups of type $A_{n}$, $B_{n}$ and $D_{n}$ is multiplicity
free. This is a useful, although mysterious, classical fact,
which is easily proved once we know it for the symmetric groups,
as $B_n$ is a wreath product and $D_n$ is a subgroup of index 2 of $B_n$.
This generalizes to the following also classical fact :

\begin{fact} If $W$ is a finite irreducible Coxeter group, it admits
a maximal parabolic subgroup $\WH$ such that the restriction
to $\WH$ of any irreducible representation of $W$ is multiplicity
free, except if $W$ has type $E_8$ or $H_4$.
\end{fact}

In case $W$ has type $E_8$ or $H_4$, there are a number
of irreducible representation whose restriction to maximal
parabolic subgroups of types $E_7$ and $H_3$ have irreducible
components with multiplicity 2. This is the worst case
scenario, so the above observation can be refined :

\begin{fact} If $W$ is a finite irreducible Coxeter group, it admits
a maximal parabolic subgroup $\WH$ such that the restriction
to $\WH$ of any irreducible representation of $W$ contains multiplicities of order at most 2,
and is even multiplicity
free, except if $W$ has type $E_8$ or $H_4$.
\end{fact}

A first goal of this note is to prove a analogous result for
the more general setting of irreducible (finite) complex pseudo-reflection groups. Recall
that such groups belong to either a finite set of 34 exceptions
or to an infinite family with three integer parameters $G(de,e,r)$.
In this family, two families can be thought of as generalisations
of Coxeter groups. The first one is when $e = 1$ : the group
$G(d,1,r)$ is a wreath product that generalizes $B_n = G(2,1,n)$.
The second one is for $d = 1$ : the groups $G(e,e,r)$ generalize both
$D_n = G(2,2,n)$ and the dihedral groups $I_2(e) = G(e,e,2)$.
Another noticeable fact, which generalizes the relation between $D_n$
and $B_n$, is that $G(de,e,r)$ is a normal subgroup of index $e$ of $G(de,1,r)$
with cyclic quotient.

It follows that the classical case-by-case approaches to the
representation of complex reflection groups and their cyclotomic
Hecke algebras usually starts with the
wreath products $G(d,1,r)$ and then use an avatar of Clifford theory
to deal with the more general groups $G(de,e,r)$ (see e.g \cite{RR,MARINMICHEL}). This approach
is however not always satisfactory. To understand this,
we can remember that many results about Coxeter groups and root systems are
simpler to prove and/or state for groups of type ADE, which have a single
conjugacy class of reflections, and then
extended or generalized to the other cases, including types $B$. The
analagous approach to complex reflection groups would be
to deal first with the groups which have a single class of reflections,
and these groups are the groups $G(e,e,r)$.
In particular, in order to generalize the above facts the crucial
case concerns the groups $G(e,e,r)$. 

\subsection{Main results}

To make the next statements precise, we need to recall some
terminology about finite complex (pseudo-)reflection groups. Let $V$
be a finite-dimensional complex vector space. A pseudo-reflection
of $V$ is an element $s \in \GL(V)$ of finite order such that
$\Ker(s-1)$ is an hyperplane of $V$. A finite subgroup $W$
of $\GL(V)$ is called a \emph{reflection group} if it is generated
by pseudo-reflections. It is called irreducible if its action
on $V$ is irreducible. A \emph{reflection subgroup} of $W$ is
a subgroup of $W$ generated by pseudo-reflections. A \emph{maximal
parabolic subgroup} of $W$ is the subgroup $W_v$ of the elements of $W$
which stabilize some given $v \in V \setminus \{ 0 \}$. It is a classical
result due to Steinberg that $W_v$ is a reflection subgroup
of $W$, generated by the pseudo-reflections of $W$ which stabilize $v$.

Recall that a matrix is called monomial
if it admits exactly one non-zero entry in each row and in each column.
Let $d,e,r \geq 1$ be integers. The group $G(de,e,r+1)$ is the subgroup of $\GL_{r+1}(\C)$ of
the monomial matrices with non-zero entries in $\mu_{de}$ such
that the product of these entries lies in $\mu_e$. The maximal
parabolic subgroup of elements leaving the $(r+1)$-th coordinate
unchanged can obviously be identified with the reflection group
$G(de,e,r)$. We refer to \cite{ARIKI,ARIKIKOIKE} for a general
account on these groups. It is known and easily checked that they are
irreducible, provided $de \neq 1$ and $(d,e) \neq (1,2)$.

We will then prove the following

\begin{theointr} The induction table between the group
$G(de,e,r+1)$ and their maximal parabolic subgroups
of type $G(de,e,r)$ contains multiplicity at most 2.
\end{theointr}

Moreover, these multiplicities appear in a systematic way that
we describe. A consequence is the following.

\begin{theointr} Any irreducible complex reflection group $W$
admits a maximal parabolic subgroup $W_v$ such that the restriction
to $W_v$ of any irreducible representation of $W$ contains
multiplicities of order at most 2, except if $W$ has type $G_{22}$, $G_{27}$.
In these cases, $W$ admits a maximal parabolic subgroup for which
the multiplicities have order at most 3.
\end{theointr}

To deduce this result from the former one, we only need to check
it for the exceptional complex reflection groups which are not
Coxeter groups. We used computer means, namely the GAP package CHEVIE.
In Table \ref{tableexc} we list all these complex reflection groups,
giving a set of generators
for a maximal parabolic subgroup satisfying our conditions,
where the names of the generators follow the conventions of the tables
in \cite{BMR}.
In all cases there exists such a subgroup which can
be generated by a subset of the usual generators, which makes
things easier to describe.
In the case of $G_{22}$ and $G_{27}$ we checked that no other maximal
parabolic subgroup behaves in a nicer way.

In the case of the groups $G(de,e,r+1)$, and in order
to be more specific about which representations of $G(de,e,r)$ occur
with multiplicity 2 in the restriction of an irreducible
representation of $G(de,e,r+1)$, we get several other results to understand
the `` square of inclusions ''
$$
\xymatrix{
G(de,1,r) \ar[r] & G(de,1,r+1) \\
G(de,e,r) \ar[u] \ar[r] & G(de,e,r+1) \ar[u]}
$$
in representation-theoretic terms. These technical results are listed
and proved in section 3.

\subsection{Representations and necklaces}

In order to prove these results for the groups $G(de,e,r)$,
we translate the questions in terms of combinatorial
data, which are called \emph{necklaces}. In general,
a necklace is a function from a group $\Gamma$,
usually assumed to be cyclic, to some set of ornaments, that
can be called pearls or colours -- and is considered
modulo the $\Gamma$-action. It is now well-known that representations
of $G(de,e,r)$ are naturally indexed by such objects (see e.g. \cite{HR}).
It turns out that understanding the branching problem involves the
following strange problem : what happens when one changes (exactly) one pearl
in a necklace ?

We did not find occurences of this problem in the
literature. Because we found it interesting in its own right,
we tried to solve it in some generality. As a consequence,
the reader interested in the proofs of the statements
in section 3 may prefer to read before that the sections 4, 6 and 7,
which deal with necklaces in general,
as a whole. Section 2 deals with a simple
general result that we use in section 3 but which does not
involve necklaces.  Section 5 contains preliminary lemmas about
cyclic groups which are used in section 6 and 7.

\bigskip
\noindent {\bf Acknowledgments.} This paper benefited from
discussions and common work about
complex reflection groups with Jean Michel, who also found
a first proof of proposition \ref{propbrisure} in case $G$ is cyclic.

\section{A general symmetry breaking result}

The aim of this section is to prove the following result, which we may have independant
interest, and that we
view as a combinatorial symmetry breaking result. For a group $G$
acting on a set $E$, and $x \in E$, we let $G_x \subset G$ denote the
stabilizer of $x$.
\begin{prop} \label{propbrisure}
Let $X$ be a set, $G$ a
group which does not contain a free product, and $X^G$ denote the set of functions
from $G$ to $X$, endowed with the natural action of $G$.
Let $\alpha,\beta \in X^G$ such that there exists a
unique $g_0 \in G$ with $\alpha(g_0) \neq \beta(g_0)$.
Then $G_{\alpha} \neq \{ 1 \} \Rightarrow G_{\beta} = \{ 1 \}$.
\end{prop}

First note that the condition on $G$ is optimal. Indeed, we can
construct a counterexample whenever $G$
contains a non-trivial free product $F = A * B$. Let $X = G$
endowed by the left-multiplication $G$-action,
and pick $a \neq b$ in $X=G$.
We let $\alpha(e) = a$, $\beta(e) = b$. Assume $w \in F \setminus \{ e \}$.
In the decomposition of $w$ in the free product $A * B$, if
the rightmost syllabon lies in $A$ we let $\alpha(w) = \beta(w) = a$,
and otherwise $\alpha(w) = \beta(w) = b$. Finally, let e.g.
$\alpha(w) = \beta(w) = a$ for all $w \in G \setminus F$. It is
easily checked that $G_{\alpha} \supset A$ and $G_{\beta} \supset B$,
although there exists a unique $x = g_0 = e \in X$
such that $\alpha(x) \neq \beta(x)$.

We remark that the condition on $G$ is closely related to the condition
of not having a free subgroup of rank 2, but is not equivalent
to it. Indeed, recall that a nontrivial free product $A * B$
contains the commutator subgroup $(A,B)$ which is free on the set $\{ aba^{-1}b^{-1}
\ | \ a \in A\setminus \{e \}, b \in B \setminus \{ e \} \}$
(see e.g. \cite{ROB} \S 6.2 exercise 7). In particular, any nontrivial
free product contains a free subgroup of rank 2, except for the
infinite dihedral group $\Z/(2) * \Z/(2) = \Z/(2) \ltimes \Z$ which does not. For groups
satisfying the Tits alternative, this condition can thus be translated
as $G$ being virtually solvable but not containing any infinite dihedral
group.

To prove this result in its full generality, we first need a criterium for a group to be a free product,
which we did not find in the literature and which has a somewhat
different flavour than the more common pingpong lemma.

\begin{prop} Let $G$ be a group acting freely on a set $E$, let
$A, B$ be two subgroups generating $G$, $K,L \subset E$ such that
$K \cap L = \{ g_0 \}$, $K^*= K \setminus \{ g_0 \}$, $L^* = L \setminus
\{ g_0 \}$. If $AK \subset K$, $BL \subset L$, $A L^* \subset L^*$
and $A K^* \subset K^*$ then $A \cap B = \{ e \}$ and $G = A * B$.
\end{prop}
\begin{proof}
We let $A^* = A \setminus \{ e \}$, $B^* = B \setminus \{ e \}$. We have
$A^* g_0 \subset K^*$, $B^* g_0 \subset L^*$ hence $A^* \cap B^* = 
\emptyset$ that is $A \cap B = \{ e \}$.

Let $\varphi : A * B \onto G$ be the natural morphism. We want to
show that $\varphi$ is injective. Let $w = A * B$. If $w$
can be written as $X_s Y_s X_{s-1} Y_{s-1} \dots X_1 Y_1$ for
some $s \geq 1$, with $X_i \in A^*$ for $i<s$, $Y_i \in B^*$
for $i \leq s$, and $X_s \in A$, we let $l(w) = s$ ; similarly
$l(w) =s$ is $w = Y_s X_s Y_{s-1} X_{s-1} \dots Y_1 X_1$ for
some $s \geq 1$, with $X_i \in A^*$ for $i\leq s$, $Y_i \in B^*$
for $i < s$, and $Y_s \in B$ ; finally $l(e) = + \infty$.

Assume by contradiction that $\Ker \varphi \neq \{ e \}$. Then
$s = \min l(\Ker \varphi) \in \Z_{>0}$ is reached for some
$w_0 \in \Ker \varphi \setminus \{ e \}$. Up to interchanging
$A$ and $B$ we may assume $w_0 = X_s Y_s \dots X_1 Y_1$ with
with $X_i \in A^*$ for $i<s$, $Y_i \in B^*$
for $i \leq s$, and $X_s \in A$.

We let $\overline{w_0} = \varphi(w_0)$, $x_i = \varphi(X_i)$,
$y_i = \varphi(Y_i)$. Note that $s \geq 2$ otherwise $y_1 = x_1^{-1} \in
A \cap B = \{ e \}$ contradicting $y_1 \neq e$. We have
$x_s y_s \dots x_1 y_1 . g_0$ hence $y_s x_{s-1} y_{s-1} \dots x_1 y_1.g_0
= x_s^{-1}.g_0$. If $x_s = e$, it follows that
$y_s x_{s-1}y_{s-1} \dots x_1 y_1 = e$ hence
$x_{s-1}y_{s-1} \dots x_1 (y_1 y_s^{-1}) = e$ ; then
$w_1 = X_{s-1} Y_{s-1} \dots X_1 (Y_1 Y_s^{-1}) \in \Ker \varphi
\setminus \{ e \}$, and $l(w_1) = s-1 < s$, a contradiction.

We thus have $x_s \in A^*$, hence $x_s^{-1}.g_0 \in K^*$. We prove
by induction that the element $y_r x_{r-1} y_{r-1} \dots x_1 y_1.g_0$
lies in $L^*$
for $1 \leq r \leq s$. The case $r = 1$ is a consequence of $B^* g_0
\subset B^* L \subset L^*$. Assuming the assertion proved for $r$,
if $r+1 \leq s$ we let $u = y_{r-1} x_{r-2}y_{r-2} \dots x_1 y_1 .g_0 \in
L^*$. Since $A^* L^* \subset L^*$ we have $x_{r-1}.u \in L^*$, and
$y_r x_{r-1}.u \in L$ since $BL \subset L$. If $y_r x_{r-1}.u = g.0$
we would have $y_r x_{r-1}y_{r-1}\dots x_1 y_1 = e$ hence
$x_{r-1} y_{r-1} \dots x_1 (y_1 y_r^{-1})$ contradicting once again
the minimality of $s$. It follows that $y_r x_{r-1}.u \in L^*$ and
we conclude by induction.

In particular, for $r = s$ we proved that $x_{s}^{-1} .g_0 \in K^* \cap
L^* = \emptyset$, a contradiction. It follows that $\varphi$ is injective
and $G = A * B$.

\end{proof}

We can now prove the main result of this section.

\begin{proof}
By contradiction we assume $G_{\alpha} \neq \{ 1 \}$
and $G_{\beta} \neq \{ 1 \}$. Let $a = \alpha(g_0)$,
$b = \beta(g_0)$,  $K = \alpha^{-1}(\{ a \})$, $K^* = \beta^{-1}(\{ a \})
= K \setminus \{ g_0 \}$ and similarly $L = \beta^{-1}(\{ b \})$, $L^* = \alpha^{-1}(\{ b \})
= L \setminus \{ g_0 \}$. It is clear that $G_{\alpha} K \subset K$,
$G_{\beta} L \subset L$ and $K \cap L = \{ g_0 \}$.

We claim that $G_{\beta} K^* \subset K^*$. Indeed, let $u \neq g_0$ in $K$
and $g \in G_{\beta}$. Then $\alpha(u) = \beta(u)$ because $u \neq g_0$
and $\beta(u) = \beta(gu)$ because $g \in G_{\beta}$. Now, if
$gu \not\in K$ then $gu \neq g_0$ hence $\beta(gu) = \alpha(gu)$
whence $\alpha(gu) = \alpha(u)$ and $gu \in K$, a contradiction.
It follows that $gu \in K$. Moreover, $gu = g_0$ would imply that
$u = g^{-1} g_0$ would satisfy both $\beta(u) = \beta(g^{-1} g_0) = \beta(g_0)$,
since $g \in G_{\beta}$, and $\beta(u) = \alpha(u)$ because $u \neq g_0$,
hence $\beta(u) = \alpha(g_0)$, contradicting $\alpha(g_0) \neq \beta(g_0)$.
The claim follows.

In the same way, $G_{\alpha} L^* \subset L^*$. By the criterium above
it follows that the subgroup of $G$ generated by $G_{\alpha}$
and $G_{\beta}$ is the free product of both, contradicting the
assumption on $G$.
\end{proof}

This result will be applied here only for a commutative group $G$, in
which case the proof does not need the criterium above. Indeed,
taking $g \in G_{\alpha} \setminus \{ 1 \}$
and $g' \in G_{\beta} \setminus \{ 1 \}$, we have $g'g g_0 \in K^*$
and $g g' g_0 \in L^*$, a contradiction since $gg'=g'g$ and $K^* \cap L^*
= \emptyset$.

\section{Representations and necklaces}

For $m \geq 2$, $d,e \geq 1$ such that $m = de$ we let
$$
G_m = \bigsqcup_{r=0}^{\infty} \Irr G(m,1,r), \ \ \ \ 
G_{d,e} = \bigsqcup_{r=0}^{\infty} \Irr G(de,e,r)
$$
and let $L : G_m \to \N = \Z_{\geq 0}$ be the map
$\rho \in \Irr G(m,1,r)  \mapsto r$. Similarly and by
abuse of notation we also denote $L : G_{d,e} \to \N$
the map $\rho \mapsto r$.

Let $\Gamma = \Z/m\Z$, $\Gamma' = d \Gamma$, $E_m = \{ X \to Y \}$
where $X = \Gamma$, viewed as a simply transitive $\Gamma$-set,
and $Y$ is the set of all partitions. There is a natural
action (on the left) of $\Gamma$ on $E_m$, given by $(\gamma.c)(x)
= c(\gamma^{-1}.x)$. For $c \in E_m$ we let
$\Aut(c) \subset \Gamma$ denote the stabilizer of $c$ in $\Gamma$.
There is a natural coding of $G_m$ by $m$-tuples of partitions
(see e.g. \cite{Zelevinsky}),
hence a natural bijective map $\Phi : G_m \to E_m$.
We have a natural map $L : E_m \to
\N$ defined by $L(c) = \sum_{x \in X} |c(x)|$ where
$|\la|$ denotes the size of the partition $\la$. This abuse of notation
is justified by $\Phi \circ L = L$.


Let $t$ be a generator of $G(m,1,1) \simeq \Z/m\Z$.
There are natural inclusions $G(m,1,r) \subset G(m,1,r+1)$ hence
$t \in G(m,1,r)$ for all $r \geq 1$. We let $t' = t^d$. The
image of $t'$ generates the cyclic quotient $G(de,1,r)/G(de,e,r) \simeq
\Z/e\Z$. Let $\zeta \in \C^{\times}$ be primitive $e$-th
root of unity. There exists a well-defined character $\eps :
G(de,1,r) \to \C^{\times}$ with kernel $G(de,e,r)$ such that $\eps(t') = \zeta$.
It is a classical fact (see e.g. \cite{HR}) that $\zeta$ can be chosen such that,
for all $\rho \in G_m$, we have $\Phi(\rho \otimes \eps) = \bar{d}.\Phi(\rho)$.
Let $r = L(\rho)$. Clifford theory says that, for $\rho_1,\rho_2 \in G_m$
with $L(\rho_i) = r$, the restrictions to $G(de,e,r)$ of
$\rho_1$ and $\rho_2$ are isomorphic iff $\rho_2 \simeq \rho_1 \otimes
\eps^n$ for some $n \in \N$, that is if $\Phi(\rho_1)$ and
$\Phi(\rho_2)$ lies in the same $\Gamma'$-orbit. On the other
hand, if $\rho \in G_{d,e}$ there exists $\tilde{\rho} \in G_m$
such that $\rho$ embeds in the restriction of $\tilde{\rho}$, and
two such $\tilde{\rho}$ are conjugated by some power of $t'$ ; in
particular they have the same restriction to $G(de,e,r)$ and,
denoting $\bar{x}$ the image of $x \in E_m$ in $E_m/\Gamma'$,
it follows that there exists a well-defined map $\overline{\Phi} :
G_{d,e} \to E_m/\Gamma'$ which sends $\rho$ to $\overline{\Phi(\tilde{\rho})}$.
Moreover, the preimage of $\bar{c} \in E_m/\Gamma'$ by $\overline{\Phi}$
has $\# \{ \tilde{\rho} \otimes \eps^n \ | \ n \in \N \}$
elements,
that is $\# \Aut_{\Gamma'}(c)$ elements, where $\Aut_{\Gamma'}(c)
= \Aut(c) \cap \Gamma'$.


The set $Y$ of partitions $\la = (\la_1 \geq \la_2 \geq \dots)$ is naturally endowed with a size function $\la \mapsto
|\la|= \sum \la_i$ and of the following usual binary relations : (non-)equality,
a total (lexicographic) ordering  $\leq$, and the relation $\la \nearrow \mu$, common
in the combinatorial representation theory of the symmetric groups, which
means $\forall i \ \la_i \leq \mu_i$ and $|\mu| = |\la|+1$. In particular
$\la \nearrow \mu$ implies $\la < \mu$.

The set $E_m$ inherits from these
the following binary relations :
\begin{itemize}
\item $\alpha \perp \beta$ if $\exists ! x \in X \ \ \alpha(x) \neq \beta(x)$;
\item $\alpha < \beta$ if $\forall x \ \alpha(x) \leq \beta(x)$ and $\alpha \perp \beta$;
\item $\alpha \nearrow \beta$ if $\alpha \perp \beta$ and
$\exists x \in X \ \ \alpha(x) \nearrow \beta(x)$.
\end{itemize}
Note that these relations are listed from the coarser to the thiner, that the
first one is symmetric and that $<$ is not a strict ordering.
In general, for an arbitrary $\Gamma$-set $X$ and $E = \{ X \to Y \}$,
with $\Gamma$ acting freely on $X$, the set of necklaces $E/\Gamma$ will be said to have \emph{ordered pearls} if
$Y$ is given a total ordering, and \emph{rough pearls} otherwise. The corresponding
combinatorics is dealt with in section \ref{sectordered} for the former case,
in section \ref{sectrough} for the latter. The relation $\perp$ is always
available, while the relation $<$ needs ordered pearls.

Let $W_r = G(de,e,r)$, $\WT_r = G(de,e,r)$.
A combinatorial
description of the branching rule for the pair $(\WT_r,\WT_{r+1})$ is
\begin{equation}
\Res_{\WT_r}^{\WT_{r+1}} \rho = \bigoplus_{\Phi(\psi) \nearrow \Phi(\rho)} \psi
= \bigoplus_{\alpha \nearrow \Phi(\rho)} \Phi^{-1}(\alpha)
\end{equation}
(see \cite[p. 104]{Zelevinsky}.) By Clifford theory and the discussion above,
a combinatorial description of the
branching rule for the pairs $(W_r,\WT_r)$ is given by
\begin{equation}
\Res_{W_r}^{\WT_{r}} \rho = \bigoplus_{\varphi \in \overline{\Phi}^{-1}(\overline{\Phi(\rho)})} \varphi
\end{equation}
We say that a representation $\rho \in \Irr(W_{r})$ \emph{extends to}
$\widetilde{W}_{r}$ if there exists $\tilde{\rho} \in \Irr
(\widetilde{W}_{r})$ such that $\rho = \Res_{W_{r}}^{\widetilde{W}_{r}}
\tilde{\rho}$.

\def\into{\mbox{ is an irreducible component of }}

\begin{prop} If $\rho_1 \in \Irr(W_{r+1})$ does not extend to $\WT_{r+1}$ then
any $\rho_2 \in \Irr(\WH)$ such that $(\Res_{W_r}^{W_{r+1}} \rho_1 | \rho_2)
\neq 0$ extends to $\WTH$. Conversely, if for $\rho_1 \in \Irr(W_{r+1})$ there
exists some $\rho_2 \in \Irr(\WH)$ not extending to $\WTH$ such
that $\rho_2 \into \Res_{\WH}^{W_{r+1}} \rho_1$, then $\rho_1$ extends
to $\WT_{r+1}$.
\end{prop}
\begin{proof} Let $\tilde{\rho}_1 \in \Irr(\WT_{r+1})$ such that
$\rho_1 \into \Res_{W_{r+1}}^{\WT_{r+1}}\tilde{\rho}_1$ and $c_1 = \Phi(\tilde{\rho}_1)$.
If $\rho_2 \in
\Irr(\WH)$ is such that $(\Res_{\WH}^{W_{r+1}} \rho_1 | \rho_2) \neq 0$,
with $\Phi(\rho_2) = \bar{c}_2$ for some $c_2 \in E_m$,
then
$$0 \neq (\Res_{\WH}^{\WT_{r+1}} \tilde{\rho}_1 | \rho_2) =
(\Res_{\WH}^{\WTH} \Res_{\WTH}^{\WT_{r+1}} \tilde{\rho}_1 | \rho_2)
= (\Res_{\WTH}^{\WT_{r+1}} \tilde{\rho}_1 | \Ind_{\WH}^{\WTH} \rho_2)
$$
meaning that we can choose $c_2 \in E_m$ such that $c_2 \nearrow c_1$.
The first assumption states $\Aut_{\Gamma'}(c_1)
\neq 1$ hence $\Aut(c_1) \neq 1$. But then
$c_2 \nearrow c_1$ implies $\Aut(c_2) = 1$ by proposition \ref{propbrisure}
whence
$\Aut_{\Gamma'}(c_2) = 1$ and
$\rho_2 = \Res_{\WH}^{\WTH} \Phi^{-1}(c_2)$.

The converse assumption states $\Aut_{\Gamma'}(c_2) \neq 1$ hence
$\Aut(c_1) \neq 1$. But then
$c_2 \nearrow c_1$ hence $c_1 \perp c_2$ and
$\Aut(c_1) = 1$ by proposition \ref{propbrisure}
whence
$\Aut_{\Gamma'}(c_1) = 1$ and
$\rho_1 = \Res_{\WH}^{\WTH} \Phi^{-1}(c_1)$.
\end{proof}

Note that $\Res_{\WH}^{W_{r+1}} \Res_{W_{r+1}}^{\WT_{r+1}} = \Res_{\WH}^{\WT_{r+1}}
= \Res_{\WH}^{\WTH} \Res_{\WTH}^{\WT_{r+1}}$ hence (1) and (2) imply
\begin{equation}
\Res_{\WH}^{W_{r+1}} \Res_{W_{r+1}}^{\WT_{r+1}} \tilde{\rho} = 
\bigoplus_{\alpha \nearrow \Phi(\tilde{\rho})} \Res_{\WH}^{\WTH}
\Phi^{-1}(\alpha) = \bigoplus_{\alpha \nearrow \Phi(\rho)}
\bigoplus_{\varphi \in \overline{\Phi}^{-1}(\overline{\alpha})} \varphi
\end{equation}
Let $\tilde{\rho_1} \in \Irr(\WT_{r+1})$. Then, for all $\rho_2 \in \Irr(\WH)$,
$$
(\Res_{\WH}^{\WT_{r+1}} \tilde{\rho}_1 | \rho_2) = \sum_{\alpha \nearrow
\Phi(\tilde{\rho}_1)} \sum_{\varphi \in \overline{\Phi}^{-1}(\overline{\alpha})}
(\varphi|\rho_2) =
\sum_{\alpha \nearrow
\Phi(\tilde{\rho}_1)} (\alpha \in \overline{\Phi}(\rho_2)) = \# \{ \alpha
\in \overline{\Phi}(\rho_2) | \alpha \nearrow \Phi(\tilde{\rho}_1) \}
$$
In particular, if $\rho_1 = \Res_{W_{r+1}}^{\WT_{r+1}} \tilde{\rho}_1$ is irreducible,
we have
$$
(\Res_{\WH}^{W_{r+1}} \tilde{\rho}_1 |\rho_2) = \# \{ \alpha \in \overline{\Phi}
(\rho_2) \ | \ \alpha \nearrow \Phi(\tilde{\rho}_1) \}.
$$ 
Otherwise, by the previous proposition we know that $\rho_2 = 
\Res_{\WH}^{\WTH} \tilde{\rho}_2$ for some
$\tilde{\rho}_2 \in \Phi(c_2)$. Let $\tilde{\rho}_1 \in \Irr(\WT_{r+1})$
such that $\rho_1 \into \Res_{W_{r+1}}^{\WT_{r+1}} \tilde{\rho}_1$ and
$c_1 = \Phi(\tilde{\rho}_1)$. We have
$$
\Res_{W_{r+1}}^{\WT_{r+1}} \tilde{\rho}_1 = \rho_1^{(1)} + \dots + \rho_1^{(s)}
$$
with $\rho_1^{(1)} = \rho_1$, $s = \# \Aut_{\Gamma'}(c_1)$.
For all $1 \leq i,j \leq s$, there exists $u$ such that
$\rho_1^{(j)} =\rho_1^{(i)} \circ \Ad(t^u)$. On the other hand,
$\rho_2 = \rho_2 \circ \Ad(t)$ hence
$$
(\Res_{\WH}^{W_{r+1}} \rho_1^{(j)} | \rho_2 ) = ( \Res_{\WH}^{W_{r+1}}
\rho_1^{(i)} \circ \Ad(t^u) | \rho_2) =
( \Res_{\WH}^{W_{r+1}}
\rho_1^{(i)}  | \rho_2\circ \Ad(t^{-u})) = 
( \Res_{\WH}^{W_{r+1}}
\rho_1^{(i)}  | \rho_2).$$
It follows that $(\Res_{\WH}^{\WT_{r+1}} \tilde{\rho}_1 | \rho_2) =
s (\Res_{\WH}^{W_{r+1}} \rho_1 | \rho_2 )$. We thus proved the
following.

\begin{prop} \label{propmultinum}If $\rho_1 \in \Irr(W_{r+1})$, $\rho_2 \in \Irr(\WH)$ with $\overline{\Phi}
(\rho_1) = \overline{c_1}$, then
$$
\left( \Res_{\WH}^{W_{r+1}} \rho_1 | \rho_2 \right) = \frac{
\# \{ \alpha \in \overline{\Phi}(\rho_2) \ | \ \alpha \nearrow c_1 \}}
{\# \Aut_{\Gamma'}(c_1)}
$$
\end{prop}

We are now ready to prove the main theorem, using combinatorial
results to be proved in the sequel.

\begin{theo} \label{theomult2} Let $\rho_1 \in \Irr(W_{r+1})$, $\rho_2 \in \Irr(\WH)$. Then
$(\Res_{\WH}^{W_{r+1}} \rho_1 | \rho_2) \leq 2$. Moreover,
if $(\Res_{\WH}^{W_{r+1}} \rho_1 | \rho_2) = 2$ then
$\rho_1$ extends to $\WT_{r+1}$.
\end{theo}
\begin{proof} Let $c_1 \in E_m$ chosen such that $\overline{c_1}=
\overline{\Phi}(\rho_1)$. If $\Aut_{\Gamma'}(c_1) = 1$,
we have to prove $\# \{ \alpha \in \overline{\Phi}(\rho_2) \ | \ \alpha
\nearrow c_1 \} \leq 2$, which is a consequence of proposition \ref{propmultmax2}.
We
thus assume $\# \Aut_{\Gamma'}(c_1) \neq 1$. Let
$
\{\alpha_1,\dots,\alpha_r \} = \{ \alpha \in \overline{\Phi}(\rho_2) \ | \ 
\alpha \nearrow c_1 \}$. Since $\overline{\Phi}(\rho_2)$ is a $\Gamma'$-orbit
we have well-defined and distincts $\gamma_i \in \Gamma' \setminus \{ 1 \}$
for $2 \leq i \leq r$ such that $\alpha_i = \gamma_i. \alpha_1$.
By lemma \ref{lemneckrestsym} we have $\gamma_i \in \Aut(c_1)$ hence
$\gamma_i \in \Aut_{\Gamma'}(c_1) = \Aut(c_1) \cap \Gamma'$. Thus
$\# \Aut_{\Gamma'}(c_1) \geq \# \{ \alpha \in \overline{\Phi}(\rho_2) \ | \ 
\alpha \nearrow c_1 \}$ and $(\Res_{\WH}^{W_{r+1}} \rho_1 | \rho_2) \leq 1$.
\end{proof} 

\begin{prop} \label{multicompo} Let $\rho \in \Irr(W_{r+1})$.
If $\rho_1,\dots, \rho_s \in \Irr(\WH)$ do not extend to $\WTH$
and satisfy that, for all $i$,
$\rho_i \into \Res_{\WH}^{W_{r+1}} \rho$, then $\forall i \  
(\Res_{\WH}^{W_{r+1}} \rho | \rho_i) = 1$ and
there exists
$\tilde{\rho}_0 \in \Irr(\WTH)$ such that each $\rho_i \into
\Res_{\WH}^{\WTH} \tilde{\rho}_0$. In particular, for all
$i,j$ there exists $g \in \WTH$ such that $\rho_j \simeq \rho_i \circ \Ad(g)$.
\end{prop}
\begin{proof}
Let $c \in E_m$ such that $\bar{c} = \overline{\Phi}(\rho)$. The statement
is void if $s \leq 1$, hence we assume $s\geq 2$. For $i \in [1,s]$,
by proposition \ref{propmultinum}, the fact that $(\Res_{\WH}^{W_{r+1}} \rho|\rho_i) \neq 0$
implies the existence of $\alpha_i \in \overline{\Phi}(\rho_i)$
such that $\alpha_i \nearrow c$. Moreover, we have $\Aut_{\Gamma'}(\alpha_i) \neq 1$
hence $\Aut(\alpha_i) \neq 1$. Then proposition \ref{proponechild} implies
$\alpha_1 = \dots = \alpha_s$, which proves the existence of
$\tilde{\rho}_0 = \Phi^{-1}(\alpha_1)$. Moreover, $\Aut(\alpha_i) \neq 1$
implies $\Aut(c) = 1$ by proposition \ref{propbrisure}. Then
proposition \ref{propmultinum} states $(\Res_{\WH}^{W_{r+1}} \rho|\rho_i ) =
\{ \alpha \in \overline{\Phi}(\rho_i) \ | \ \alpha \nearrow c \}$,
hence $(\Res_{\WH}^{W_{r+1}} \rho | \rho_i) = 1$ again by proposition \ref{proponechild}.
The final assertion is an immediate
consequence of Clifford theory.
\end{proof}

We may wonder for which $\rho \in \Irr(W_{r+1})$ there exists several $\rho_1, \dots, \rho_s
\in \Irr(\WH)$ such that $(\Res_{\WH}^{W_{r+1}} \rho | \rho_i) = 2$. By theorem
\ref{theomult2}, $\rho$ extends to some $\tilde{\rho} \in \Irr(\WT_{r+1})$
and by proposition \ref{multicompo} each $\rho_i$ extends to some
$\tilde{\rho}_i \in \Irr(\WTH)$. Let $\alpha_i \in \Phi(\tilde{\rho}_i)$,
$c = \Phi(\tilde{\rho})$. We have $\Aut_{\Gamma'}(\alpha_i) = \Aut_{\Gamma'}
(c) = 1$, $\alpha_i \nearrow c$, and there exists $\gamma_i \in \Gamma'
\setminus \{ e \}$ such that $\gamma_i . \alpha_i \nearrow c$ with $r \geq 2$.
Then the possible shapes of $c$ are implicitely given by proposition
\ref{propmulticompo2}.

They are most easily described when $d = 1$, that is $\Gamma = \Gamma'$
and $W = G(e,e,r+1)$. Choose some nontrivial subgroup $\Gamma_0 = m \Gamma$
of $\Gamma$, and define $c : \Gamma \to Y$ as follows. Subdivide
$\Gamma = \Gamma_0 \sqcup (x_1+\Gamma_0) \sqcup \dots \sqcup (x_{m-1} +
\Gamma_0)$ in $\Gamma_0$-cosets, pick one $\la_i \in Y$ for
each $i \in [1,m-1]$ and put $c(x) = \la_i$ for $x \in x_i + \Gamma_0$;
choose $\la_0,\mu_0 \in Y$ such that $\mu_0 \nearrow \la_0 \in Y$
and define $c(0) = \mu_0$, $c(x) = \la_0$ for $x \in \Gamma_0 \setminus
\{ 0 \}$. If $u = e/m = | \Gamma_0 |$, it is easily checked that
$\Aut(c) = 1$ and that the restriction of $\rho = \Res_{W_{r+1}}^{\WT_{r+1}} \Phi^{-1}(c)$
to $\WH$ admits at least $\llcorner \frac{u-1}{2} \lrcorner$ components of multiplicity 2. Proposition
\ref{propmulticompo2} and its corollary state that all $\rho \in \Irr(W)$ whose restriction
to $\WH$ contains at least 2 components with multiplicity 2 are obtained
in this way, with $u \geq 5$.

\section{A basic lemma on necklaces}

Here we let $X$ be
a finite set acted upon \emph{freely} by a group $\Gamma$, $Y$
be a set containing at least 2 elements,
and $E = Y^X$ be the set of maps
from $X$ to $Y$. There is a natural action of
$\Gamma$ on $E$. We recall that for $\alpha,c \in E$
the notation $\alpha \perp c$ means
$\exists ! x \in X \ \ \alpha(x) \neq c(x)$.
and denote $\Aut(c) \subset \Gamma$ the
stabilizer of $c$ under the action of $\Gamma$.

When $\Gamma$ is cyclic of order $n$ and $\gamma \in \Gamma$,
we also introduce the following notation.
Let $v \in X$ and $u \in \Gamma.v$. If
$r \in [0,n-1]$ is defined by
$u = \gamma^r . v$, then we let
$[v,u]_{\gamma} = \{ \gamma^k .v \ | \ k \in [0,r] \}$.
The companion notations
$]v,u]_{\gamma}$, $[v,u[_{\gamma}$, $]v,u[_{\gamma}$ are
self-explaining.

The set $Y$ has to be thought of as a
set of pearls, distinguished by an ornament. In
order to clearly distinguish elements of
the $\Gamma$-set $X$ from elements of the set $Y$ we use ornamental
symbols $\spadesuit, \heartsuit,\diamondsuit ,\clubsuit$
for elements of $Y$ in the proofs.

The following technical lemma is basic for our
purposes.

\begin{lemma} \label{lemneckbasic} Assume that $\Gamma$ is cyclic
with generator $\gamma$ and acts freely on $X$.
Let $c \in E$. The following
are equivalent
\begin{itemize}
\item[(i)] $\exists \alpha,\beta \in E$ such that $\alpha \neq \beta$,
$\alpha \perp c$, $\beta \perp c$ and $\beta = \gamma.\alpha$
\item[(ii)] $\exists  \mathcal{O} \in X/\Gamma$ such that
\begin{itemize}
\item[(a)] $c$ is constant on each $P \neq \mathcal{O}$ in $X/\Gamma$.
\item[(b)] there exists $u,v \in \mathcal{O}$
such that $u \neq v$ and $c$ is constant on both $[v,u]_{\gamma}$
and its complement in $\mathcal{O}$.
\end{itemize}
\end{itemize}
Under these assumptions, we have $\gamma \in \Aut(c) \Leftrightarrow | c(\mathcal{O}) | = 1
\Leftrightarrow v = \gamma.u$. Moreover, $u$ and $v$ are
characterized in $X$ by $\alpha(u) \neq c(u)$ and $\beta(v) \neq c(v)$.
Finally we have $| \alpha(\mathcal{O}) | = | \beta(\mathcal{O})|$.
\end{lemma}
\begin{proof}
$(ii) \Rightarrow (i)$. (b) implies $| c(\mathcal{O})| \leq 2$.
We let $\clubsuit \in Y$ such that $c(\mathcal{O}) = \{ c(u),\clubsuit \}$
if $| c(\mathcal{O}) |= 2$, and let $\clubsuit$ be an arbitrarily
chosen element of $Y \setminus c(\mathcal{O})$
otherwise, using $|Y| \geq 2$. We have to define $\alpha$ and $\beta$
fulfilling $(i)$.
We define $\alpha(u) = \clubsuit$, $\alpha(x) = c(x)$ for $x \neq u$
and  $\beta(v) = \clubsuit$, $\beta(x) = c(x)$ for $x \neq v$.
We have $\alpha(u) \neq c(u)$ hence $\alpha(u) \neq \beta(u) = c(u)$
because $u \neq v$. It follows that
$\alpha \neq \beta$, and clearly $\alpha \perp c$, $\beta \perp c$.
It remains to show that $\beta = \gamma.\alpha$.

Let $x \in X$. If $\Gamma.x \neq \mathcal{O}$ then $c(\gamma.x) = c(x)$
by (a) hence $\beta(\gamma.x) = c(\gamma.x) = c(x) = \alpha(x)$ since
$x,\gamma.x \in \Gamma.x$ and $u,v \not\in \Gamma.x$. Now (b) tells
us that $\beta$ equals $c(u)$ on $[u,\gamma^{-1}.v]_{\gamma}$, $\clubsuit$
on its complement in $\mathcal{O}$, and $\alpha$ equals $c(u)$
on $[\gamma.u,v]_{\gamma}
= \gamma.[u,\gamma^{-1}.v]_{\gamma}$,
$\clubsuit$ on its complement. It follows that
$\beta(\gamma.x) = \alpha(x)$ also for $x \in \mathcal{O}$, hence
$\beta = \gamma.\alpha$.


$(i) \Rightarrow (ii)$. Since $\alpha \perp c$ and $\beta \perp c$
there exist well-defined $u,v \in X$ such that $\alpha(u) \neq c(u)$
and $\beta(v) \neq c(v)$. Let $\clubsuit = \alpha(u)$
and $e = | c^{-1}(\clubsuit)|$. If $u = v$, $\alpha \neq \beta$ implies
$\beta(u) \neq \clubsuit$ hence $|\beta^{-1}(\clubsuit) | = e$,
but $|\beta^{-1}(\clubsuit)| = |\alpha^{-1}(\clubsuit)| = e+1$, since
$\beta = \gamma.\alpha$. It follows that $u\neq v$ and $|\beta^{-1}
(\clubsuit)| = e+1$, hence $\beta(v) = \clubsuit$. Similarly,
considering $\spadesuit = c(u) \neq \clubsuit$ and
$f = |c^{-1}(\spadesuit)|$,
we get $c(u) = \beta(u) = \alpha(v) = c(v) = \spadesuit$, since
$f-1 = |\alpha^{-1}(\spadesuit)| = |\beta^{-1}(\spadesuit)|$.




We let $\mathcal{O} = \Gamma.u  \in X/\Gamma$, and $n = |\Gamma|$. Assume
by contradiction $v \not\in \mathcal{O}$. Then by induction we
have $\beta(\gamma^r.u ) = \alpha(u)$ for all $r \in [1,n]$. Indeed,
$\beta = \gamma.\alpha$ proves the case $r=1$, and also
implies $\beta(\gamma^{r+1}.u) = \alpha(\gamma^r.u)$ ; then
$r<n$ implies $\gamma^r.u \neq u$ hence
$\alpha(\gamma^r.u) = c(\gamma^r.u)$ ; $v \not\in \mathcal{O}$ implies
$c(\gamma^r.u)  = \beta(\gamma^r.u)$, which equals $\alpha(u)$
by the induction hypothesis. This yields the
contradiction $c(u) = \beta(u) = \beta(\gamma^n.u) = \alpha(u)$.
Hence $v \in \mathcal{O}$, and $u = \gamma^m.v$ for
some $m \in [1,n-1]$. For $0 \leq r \leq m-1$ we have
$\gamma^{r+1}.v \neq v$ and $\gamma^r.v \neq u$ hence
$c(\gamma^{r+1}.v) = \beta(\gamma^{r+1}.v) = \alpha(\gamma^r.v)
=c(\gamma^r.v) = c(\gamma^r.v)$ meaning that $c$ is constant on
$[v,u]_{\gamma}$. Similarly, if $v = \gamma^m .u$ for some
$m \in [1,n-1]$, then for $1 \leq r \leq m-2$ we have $\gamma^r.u \neq u$
and $\gamma^{r+1}.u \neq v$, whence
$c(\gamma^{r+1}.u) = \beta(\gamma^{r+1}.u) = \alpha(\gamma^r.u)=
c(\gamma^r.u)$ and $c$ is constant on its complement, which proves (b).
Let now $P \in X/\Gamma$ with $P \neq \mathcal{O}$. Since $u,v \not\in P$
and $\beta = \gamma.\alpha$ we have $c(\gamma.x) = \beta(\gamma.x) =
\alpha(x) = c(x)$ for all $x \in P$, which proves (a). 
The proof that $\gamma \in \Aut(c) \Leftrightarrow | c(\mathcal{O})|
=1 \Leftrightarrow v = \gamma.u$ is straightforward. Finally,
we show that $|\alpha(\mathcal{O})| = |\beta(\mathcal{O})| = 2$.
If $|c(\OO)| =1$ we have $|\alpha(\OO)| = |\beta(\OO)| = \left|
\{ \beta(v),c(v) \} \right| = 2$. We thus can assume $|c(\OO)| = 2$, which
implies $v \neq \gamma.u$. Assume by contradiction that $\beta(v)
= \heartsuit \not\in c(\OO)$. We have $\alpha(u) = \beta(v) = \heartsuit$
and, for $x \in \OO$, $x=v \Leftrightarrow \beta(x) = \heartsuit$
and $x = u \Leftrightarrow \alpha(x) = \heartsuit$. Then
$\beta(\gamma.u) = \alpha(u) = \heartsuit$ implies $\gamma.u = v$
which has been excluded. Thus $\beta(v) \in c(\OO)$ and
$|\alpha(\OO)| = |\beta(\OO)| = |c(\OO)| = 2$.

\end{proof}

\section{Preliminaries on cyclic groups}
\label{sectcyclic}

\begin{lemma} \label{lemneckclasses} Let $\Gamma$ be a cyclic group acting freely and
transitively on a finite set $X$. Let $\Gamma_1,\Gamma_2$ be subgroups of $\Gamma$
such that $\Gamma = \Gamma_1 \Gamma_2$. For all
$(P,Q) \in X/\Gamma_1 \times X/\Gamma_2$ we have $P \cap Q \neq \emptyset$.
\end{lemma}

\begin{proof} Let $n = |\Gamma|$ and $n_i = |\Gamma_i|$. Since
$\Gamma$ is cyclic we have $|\Gamma_1 \cap \Gamma_2| = \gcd(n_1,n_2)$
and the assumption $\Gamma = \Gamma_1 \Gamma_2$ means $|\Gamma|
= \lcm(n_1,n_2)$. If $(P,Q) \in X/\Gamma_1 \times X/\Gamma_2$
satisfies $P \cap Q \neq \emptyset$, then $\Gamma_1 \cap \Gamma_2$
acts freely on $P \cap Q$. Moreover, if $x,y \in P \cap Q$, we know
that there exists $\gamma_1 \in \Gamma_1$ and $\gamma_2 \in \Gamma_2$
such that $y = \gamma_1.x = \gamma_2.x$, hence $\gamma_2^{-1} \gamma_1.x = x$
and $\gamma_1 = \gamma_2 \in \Gamma_1 \cap \Gamma_2$ because $\Gamma$ acts freely on $X$.
It follows that $\Gamma_1 \cap \Gamma_2$ acts freely and transitively on
$P \cap Q$ hence $|P \cap Q| = |\Gamma_1 \cap \Gamma_2 | = \gcd(n_1,n_2)$.
Now let $P \in X/\Gamma_1$. It is the disjoint union of the $P \cap
Q$ for $Q \in X/\Gamma_2$, hence
$$
n_1 = |P| = \sum_{Q \in X/\Gamma_2} |P \cap Q| \leqslant \gcd(n_1,n_2) |X/\Gamma_2|
\leqslant \frac{\gcd(n_1,n_2) \lcm(n_1,n_2)}{n_2} = n_1
$$
and $\sum_{Q \in X/\Gamma_2} |P \cap Q| = \gcd(n_1,n_2) |X/\Gamma_2|$
hence $|P \cap Q| \neq 0$ for all $Q \in X/\Gamma_2$.
\end{proof}

We define $\Aff_n$ to be the group of bijective affine functions from $\Z/n\Z$ to
itself :
$$\Aff_n = \{ \varphi \in \mathrm{Bij}(\Z/n\Z) \ | \ 
\exists \alpha,\beta \in \Z/n\Z \ \ \forall x \in \Z/n\Z \ \ \varphi(x) = \alpha + \beta x \}.$$

We will use the following lemma.

\begin{lemma} \label{lemneckaffine} Let $n \geq 3$, $0 \leq m \leq n-2$ and $I_m = \{ \bar{0}, \bar{1} ,
\dots , \bar{m} \} \subset \Z/(n)$. Let $\varphi \in \Aff_n$ such that $\varphi(I_m) \subset I_m$.
Then :
\begin{enumerate}
\item If $1 \leq m \leq n-3$ then $\varphi = \Id$ or $\forall x \in \Z/(n)
\ \ \varphi(x) = m-x$.
\item If $m = 0$ there exists $r \in [0,n-1]$ with $\gcd(r,n) = 1$ such that
$\forall x \in \Z/(n) \ \ \varphi(x) = rx$.
\item If $m = n-2$ there exists $r \in [0,n-1]$ with $\gcd(r,n) = 1$ such that
$\forall x \in \Z/(n) \ \ \varphi(x) = r-1+rx$.
\end{enumerate}
\end{lemma}
\begin{proof}
Since $\varphi$ is injective we know that $\varphi(I_m) = I_m$. Let
$a \in [0,n-1]$ and $r \in (\Z/(n))^{\times}$ such that
$\varphi(x) = \bar{a} + r x$ for all $x \in \Z/(n)$.
If $m = 0$ then $\varphi(\bar{0}) = \bar{0}$ hence $a = 0$ and the
conclusion follows. If $m = n-2$ then $\{ - \bar{1} \}
= \Z/(n) \setminus I_m$ hence $\varphi(-\bar{1}) = - \bar{1}$
that is $\bar{a} = r-\bar{1}$ and the conclusion follows.

We thus can restrict ourselves to
assumption (1). Assume for now
that $m < n-m$. Let $\Delta : I_m\times I_m \to
\Z/(n)$ be defined by $\Delta(x,y) = x-y$. The set
$\Delta(I_m\times I_m) = \{ - \bar{m} ,\dots, -\bar{1} , \bar{0}, \bar{1},
\dots , \bar{m}\}$ has cardinality $2m +1 \leq n$. Let
$\Phi$ be the restriction of $\varphi \times \varphi$
to $I_m\times I_m$. This is a bijection of $I_m\times I_m$.
We have $| \Delta^{-1}
(y)| = m$ if and only if $y \in \{ -\bar{1},\bar{1} \}$. Since
$|(\Delta \circ \Phi)^{-1}(\bar{1})| = |\Delta^{-1}(\bar{1})|$
by bijectivity of $\Phi$ and
$(\Delta \circ \Phi)^{-1}(\bar{1}) = \Delta^{-1}(r^{-1})$
by direct calculation, it follows that
$r \in \{ -\bar{1},\bar{1} \}$. If $r = \bar{1}$ then
$\varphi(x) = \bar{a} + x$ for all $x \in \Z/(n)$.
Consider in that case the iterated maps $\varphi^j$ of $\varphi$
for $j \in \N$. These
induce bijections of $I_m$. If $a \neq 0$ there
would exist $j \in \Z_{>0}$ such that $ja > m$
and $(j-1)a \leq m$. But $\bar{ja} = \varphi^j(\bar{0}) \in I_m$
hence $ja \geq n$ and $a = ja - (j-1)a > n-m > m$
by assumption, a contradiction since $a \in [0,m]$.
It follows that $\varphi = \Id$. If $r = -\bar{1}$,
we introduce $\psi_m \in \Aff_n$
defined by $\psi_m(x) = \bar{m} -x$. Then
$\psi_m \circ \varphi \in \Aff_n$
sends $I_m$ into itself and $\psi_m \circ \varphi(x)
= \bar{m} - \bar{a} + x$, hence $\psi_m \circ \varphi
= \Id$ by the above discussion, $\bar{m} = \bar{a}$
and $a = m$ hence $\varphi = \psi_m$
since $\psi_m^2 = \Id$.

Now assume $m \geq n-m$. Let $S \in \Aff_n$
defined by $S(x) = -\bar{1} - x$ for all
$x \in \Z/(n)$, and $I'_m = \{ \overline{m+1}, \dots,
\overline{n-1} \}$. We have $\varphi(I'_m) = I'_m$.
Let $\varphi' = S \circ \varphi \circ S \in \Aff_n$.
We have $S(I_m) = I'_{n-m-2}$, $S(I'_m) = I_{n-m-2}$
hence $\varphi'(I_{n-m-2}) = I_{n-m-2}$. Moreover
$1 \leq m \leq n-3$ implies $1 \leq n-m-2 \leq n-3$.
It follows that $\varphi' \in \{ \Id, \psi_{n-m-2} \}$
since $n-m-2 < n-m \leq m < m+2 = n-(n-m-2)$
and thus $\varphi \in \{ \Id, \psi_m \}$.
\end{proof}

\section{Necklaces with rough pearls} 
\label{sectrough}

In this section we deal with the case where
the set $Y$ has no additional structure.
We recall that $\Gamma$ is cyclic and acts freely on the finite set $X$.

\begin{lemma} \label{lemneckrestsym}
Let $c \in E$ such that $\Aut(c) \neq 1$. If
there exists $\alpha \perp c$ and $\gamma \in \Gamma$
such that $\gamma.\alpha \perp c$
then $\gamma \in \Aut(c)$.
\end{lemma}
\begin{proof}
Let $\beta = \gamma.\alpha$. By assumption there exists
$\delta \in \Aut(c) \setminus \{ 1 \}$.
We assume by contradiction that $\gamma \not\in \Aut(c)$.
In particular $\gamma \neq 1$ and $\alpha \neq \beta$. Let $\Gamma_0
= < \gamma >$, $\Delta_0 = < \delta >$ and $\Gamma' = < \Gamma_0, \delta_0>$.

Lemma \ref{lemneckbasic} applied to $\Gamma_0 = < \gamma >$
defines $u,v \in X$ and $\mathcal{O} = X/\Gamma_0$. Since
$\gamma \not\in \Aut(c)$ these elements are uniquely defined.
We let $X' = \Gamma'.v$. Obviously $\mathcal{O} \subset X'$.

Since $\Gamma' = \Gamma_0 \Delta_0$ acts freely
and transitively on $X'$ we get by lemma \ref{lemneckclasses} that
$P \cap Q \neq \emptyset$ for all
$(P,Q) \in X'/\Gamma_0 \times X'/\Delta_0$. Since $\delta
\in \Aut(c)$ the map $c$ is constant on each $Q \in X'/\Delta_0$,
hence induces a map $\bar{c} : X'/\Delta_0 \to Y$.
If there were $P \neq \mathcal{O}$ in $X'/\Gamma_0$ then
$c$ would be constant on $P$, hence $\bar{c}$ and $c$
would be constant. This is a contradiction because
$c$ is not constant on $\mathcal{O} \subset X'$.
It follows that $X'/\Gamma_0 = \{ \mathcal{O} \}$
and $\Gamma' = \Gamma_0$. In particular
$\delta \in \Gamma_0$ and $\delta = \gamma^r$,
for some $r \in [2,n-2]$ since $\gamma \not\in \Delta_0 \subset \Aut(c)$.

Since
$X' = \mathcal{O}$ with $|c(\mathcal{O})| = 2$ and $c$ is constant on each $Q \in \mathcal{O}/\Delta_0$,
we know that $[v,u]_{\gamma}$ is a union of $\Delta_0$-orbits, hence
is $\delta$-stable. Let $n$ denote the order of $\gamma$. We identify $\mathcal{O}$ with $[0,n-1]$, $v$ with $0$,
$u$ with $m \in [1,n-2]$, $\gamma$ to $\bar{1} \in \Z/(n)$. Since
$\delta.v \in [v,u]$ and $r \in [2,n-1]$ we have $r \leq m$. Let
$w = \gamma^{-1}.v$, identified with $n-1$. We have $c(w) \not\in
c([v,u]_{\gamma})$ hence $\delta.w \not\in [v,u]_{\gamma}$.
On the other hand $\delta.w$ is identified with $r-1 \geq 0$,
but $r \in [2,n-2]$ implies $r-1 \geq 0$ and $r \leq m$ implies
$r-1 \leq m$. It follows that $\delta.w \in [v,u]_{\gamma}$,
a contradiction.

\end{proof}

\begin{prop} \label{propmultmax2} Let $\alpha_1,\alpha_2,\beta \in E$ such that
$\alpha_1,\alpha_2,\beta \perp c$, $\beta = \gamma_1 . \alpha_1
= \gamma_2.\alpha_2$ with $\gamma_1,\gamma_2 \in \Gamma
\setminus \Aut(c)$. Then $\alpha_1 = \alpha_2$.
\end{prop}
\begin{proof}
We assume by contradiction that $\alpha_1, \alpha_2$
and $\beta$ are all distinct. Let $\Gamma_i = < \gamma_i>$.
Since $\gamma_1,\gamma_2 \not\in \Aut(c)$, lemma
\ref{lemneckbasic} provides two special orbits
$\mathcal{O}_1, \mathcal{O}_2$ and
$u_1,u_2,v_1,v_2 \in X$ with $u_1 \neq u_2$. Since $v_1,v_2$
are characterized in $X$ by $\beta(v_i) \neq c(v_i)$ we have
$v_1=v_2=v$.

We first rule out the possibility that $\Gamma_1 = \Gamma_2$.
In that case, let $\Gamma_0 = \Gamma_1 = \Gamma_2$.
We have $\gamma_2 = \gamma_1^r$ for some $r$
prime to $n = |\Gamma_0|$, and
$\mathcal{O}_1 = \mathcal{O}_2 = \Gamma_0.v = \mathcal{O}$
can be identified with $\Z/(n)$, $v$ with $\bar{0}$, $\gamma_1$
with $\bar{1} \in \Z/(n)$. Let $\varphi : x \mapsto r x$
in $\Aff_n$. We have $[v,u_1]_{\gamma_1} = [v,u_2]_{\gamma_2}$,
which means that $\varphi$ preserves $I_m \subset \Z/(n)$
where $u_1$ is identified with $\bar{m}$ for
some $m \in [1,n-2]$ (recall that $|c(\mathcal{O})| = 2$
hence $[v,u_1]_{\gamma_1} \neq \mathcal{O}$).
Since in addition $\varphi(\bar{0}) = \bar{0}$,
lemma \ref{lemneckaffine} implies $\varphi = \Id$
meaning $\gamma_1 = \gamma_2$ and $\alpha_1=\alpha_2$,
a contradiction.

Let $\Gamma' = \Gamma_1 \Gamma_2$, $X' = \Gamma'.v$
and $c'$ the restriction of $c$ to $X'$.
We have $\mathcal{O}_1,\mathcal{O}_2 \subset X'$.
Assume that $X' \neq \mathcal{O}_1$
and $X' \neq \mathcal{O}_2$, or equivalently
$\Gamma' \neq \Gamma_1$, $\Gamma' \neq \Gamma_2$. Then there
exists $P \neq \mathcal{O}_1$ in $X'/\Gamma_1$ and
$Q \neq \mathcal{O}_2$ in $X'/\Gamma_2$. Since $P$
intersects each element of $X'/\Gamma_2$ and $Q$
intersects each element of $X'/\Gamma_1$ by lemma
\ref{lemneckclasses}, we get that $c'$ is constant
on $X' \setminus \mathcal{O}_1 \cap \mathcal{O}_2$.
Let $\heartsuit \in Y$ be the value it takes.
Since $\gamma_i \not\in
\Aut(c)$ we have $|c(\mathcal{O}_i)| = 2$.
We know that $c(X' \setminus \mathcal{O}_1) = \{ \heartsuit \}$.
On the other hand, $\mathcal{O}_1 \cap \mathcal{O}_2
\varsubsetneq \mathcal{O}_1$ otherwise $\mathcal{O}_1 \subset \mathcal{O}_2$,
in particular $\gamma_1.v \in \mathcal{O}_2$ and $\gamma_1
\in \Gamma_2$ by the freeness assumption, hence $\Gamma_1 \subset
\Gamma_2$ contradicting $\Gamma_2 \neq \Gamma'$. It follows that
$\heartsuit \in c(\mathcal{O}_1)$ and
$c(X') = c(\mathcal{O}_1) = c(\mathcal{O}_2) = \{ \heartsuit, \clubsuit \}$
for some
$\clubsuit \neq \heartsuit$.
We claim that there exists only one $x \in X'$ such
that $c(x) = \clubsuit$. By contradiction assume
otherwise. These elements belong to $\mathcal{O}_1$, hence by lemma
\ref{lemneckbasic} they belong either to $[v,u_1]_{\gamma_1}$
or to its complement $]u_1,v[_{\gamma_1}$. If there are
at least two of them, we then have some $x \in \mathcal{O}_1$
such that $c(\gamma_1.x) = c(x) = \clubsuit$. But then
$x,\gamma_1.x \in \mathcal{O}_1 \cap \mathcal{O}_2 \subset \mathcal{O}_2$
hence $\gamma_1 \in \Gamma_2$ by the freeness assumption and $\Gamma_1
\subset \Gamma_2$, a contradiction.
Let then $x$ denote the only
element in $X'$ satisfying $c(x) = \clubsuit$.
Since $c(v) = c(u_1)$ and $v \neq u_1$ we have $c(v)
= c(u_1) = \heartsuit$.
Likewise, $c(u_2) = \heartsuit$. This implies
$\gamma_1.x = v$ and $\gamma_2.x = v$,
hence $\gamma_1 = \gamma_2$, a contradiction.

It follows that $\Gamma_1 \subset \Gamma_2$
or $\Gamma_2 \subset \Gamma_1$. By symmetry
we may assume $\Gamma_1 \subset \Gamma_2$,
that is $<\gamma_1> = <\gamma_2^r>$ for some $r \in [1,n-1]$
dividing the order $n$ of $\Gamma_2$, and $r \geq 2$
since $\Gamma_1 \neq \Gamma_2$ (of course we do not necessarily have
$\gamma_1 = \gamma_2^r$).
We identify $\mathcal{O}_2$ and $\Gamma_2$ with $\Z/(n)$, $v$
with $\bar{0}$, $\gamma_2$ with $\bar{1}$. Then
$u_1 = \overline{m_1 r}$ for some $1 \leq m_1
< n/r$ since $u_1 \neq v$. Similarly
$u_2 = \overline{m_2}$ for some
$m_2 \in [1,n[$. Let $P = \bar{1} + \bar{r} \Z
\subset \mathcal{O}_2$. We have $P \cap \mathcal{O}_1
= \emptyset$ since $r \geq 2$, hence $c$ is constant on $P$. Let $c(v) =
\spadesuit$, $c(\mathcal{O}_2) = \{ \spadesuit, \clubsuit \}$.
Since $c([v,u_2]_{\gamma_2}) = \{ \spadesuit \}$,
$[v,u_2]_{\gamma_2}$ has been identified with $[0,m_2]$,
and the class of $1 \in [0,m_2]$ belongs to $P$,
we get $c(P) = \{ \spadesuit \}$.



On the other hand, $c(\mathcal{O}_1) = \{ \spadesuit,\clubsuit \}$ since
$\gamma_1 \not\in \Aut(c)$, hence there exists
$k \in [1,\frac{n}{r} [$ such that $c(n-kr) = \clubsuit$. It
follows that $n-kr \in ]u_2,v[_{\gamma_2}$ and $c(x) = \clubsuit$
for all $x \in [n-kr,n-1]$. In particular $c([n-r,n-1]) = \{ \clubsuit \}$
hence $[n-r,n-1]\cap P = \emptyset$. But $\bar{n} \not\in P$ hence
$P \cap \overline{[n-r,n]} = \emptyset$, a contradiction since
$P = 1 + \bar{r} \Z$.
\end{proof}

\section{Necklaces with ordered pearls}
\label{sectordered}

We assume here that $Y$ is endowed with a \emph{total} ordering $\leq$.
This enables one to introduce the following relation
on $E$ : we note $\alpha < \beta$ if
$\alpha \perp \beta$ and $\alpha(x) \leq \beta(x)$
for all $x \in X$. In terms of pearls, we can imagine that the
elements of $Y$ are greyscales, and that $\alpha$ is deduced
from $\beta$ by fading one pearl.
If $\alpha < \beta$ we call $\alpha$ a \emph{child}
of $\beta$. Recall that $\Aut(\alpha)$ denotes the stabilizer of $\alpha$
in $\Gamma$. We say that two childs $\alpha_1 \neq
\alpha_2$ of $c$ are \emph{twins} if there exists $\gamma \not\in \Aut(c)$
such that $\alpha_2 = \gamma.\alpha_1$. We assume again that $\Gamma$
is cyclic and acts freely on the finite set $X$. By proposition \ref{propmultmax2}
above, we know that triplets do not occur.

\subsection{At most one child admits symmetries}

\begin{prop} \label{proponechild} Let $\alpha_1,\alpha_2,c \in E$ such that
$\alpha_1,\alpha_2 < c$. If $\Aut(\alpha_1) \neq 1$ and $\Aut(\alpha_2)
\neq 1$ then $\alpha_1=\alpha_2$.
\end{prop}
\begin{proof} We argue by contradiction, assuming $\alpha_1 \neq \alpha_2$.
Let $\Gamma_1 = \Aut(\alpha_1)$, $\Gamma_2 = \Aut(\alpha_2)$ and
$x_1,x_2 \in X$ such that $\alpha_i(x_i) < c(x_i)$ for $i \in \{ 1, 2 \}$.

As a first step, we prove that this implies $x_1 \neq x_2$. 
We assume otherwise
and let $x_0 = x_1 = x_2$. An element $g \in \Gamma_1 \cap
\Gamma_2 \setminus \{ e \}$ would yield
$\alpha_1(x_0) = \alpha_1(g.x_0) = c(g.x_0) = \alpha_2(g.x_0) =
\alpha_2(x_0)$ hence $\alpha_1 = \alpha_2$, a contradiction. Let
then $g_i \in \Gamma_i \setminus \{ e \}$ for $i \in \{ 1, 2 \}$.
If $g_2 g_1.x_0 = x_0$ then $g_2 = g_1^{-1} \in \Gamma_1 \cap \Gamma_2
\setminus \{ e \}$ which has been ruled out. Thus $
\alpha_2(g_2 g_1.x_0) = c(g_2g_1.x_0)$ and
$$
c(g_2g_1.x_0) = \alpha_1(g_2g_1.x_0) = \alpha_1(
g_1g_2.x_0) = \alpha_1(g_2.x_0) = c(g_2.x_0) = \alpha_2(g_2.x_0) = \alpha_2(x_0)
$$
and also $\alpha_2(g_2 g_1.x_0)= \alpha_2(g_1.x_0)
= c(g_1.x_0) = \alpha_1(g_1.x_0) = \alpha_1(x_0)$
hence $\alpha_1 = \alpha_2$, a contradiction.

We thus proved $x_1 \neq x_2$. As a second step, we prove $\Gamma_1
\cap \Gamma_2 = \{ e \}$, by contradiction. Assume we have
$g \in \Gamma_1 \cap \Gamma_2$ with $g \neq e$, and recall $x_1 \neq x_2$.
If $x_2 \neq g.x_1$ then, on the one hand $\alpha_2(g.x_1) = c(g.x_1)=\alpha_1(g.x_1)
= \alpha_1(x_1)$, and on the other hand $\alpha_2(g.x_1) = \alpha_2(x_1)
=c(x_1)$ since $x_2 \neq x_1$, hence $c(x_1) = \alpha_1(x_1)$,
a contradiction. It follows that $x_2 = g.x_1$. This implies
$|\Gamma_1 \cap \Gamma_2| = 2$ by freeness of the $\Gamma$-action,
hence $g= g^{-1}$ and $x_1 = g.x_2$. But then follows the
following contradiction :
$$
\left\lbrace \begin{array}{lclclclcl}
c(x_1) &=& \alpha_2(x_1) &=& \alpha_2(g.x_2) &=& \alpha_2(x_2) &<& c(x_2) \\
c(x_2) &=& \alpha_1(x_2) &=& \alpha_1(g.x_1) &=& \alpha_1(x_1) &<& c(x_1). \\
\end{array} \right.
$$
As a consequence we get that, for all $(g_1,g_2) \in \Gamma_1 \times
\Gamma_2$ with $g_1,g_2 \neq e$, we have
$$
|\{ x_1,g_1.x_1,g_2.x_1,g_1g_2.x_1 \}| = 4 \mbox{ and }
|\{ x_2,g_1.x_2,g_2.x_2,g_1g_2.x_2 \}| = 4.
$$

As a third step we prove that $x_2 \not\in \Gamma_2.x_1$ and $x_1 \not\in
\Gamma_1.x_2$. By symmetry considerations it is sufficient to show that
$x_2 \not\in \Gamma_2.x_1$. We argue by contradiction, assuming
$x_2 = g_2 . x_1$ with $g_2 \in \Gamma_2$. Since $x_1 \neq x_2$
we know that $g_2 \neq e$. Moreover, this also
implies $c(x_1) = \alpha_2(x_1) = \alpha_2(g_2.x_1) = \alpha_2(x_2) < c(x_2)$
and $c(x_2) = \alpha_1(x_2) = \alpha_1(g_1.x_2)=\alpha_1(g_1g_2.x_1)$
for all $g_1 \in \Gamma_1$. By assumption we can choose $g_1 \in \Gamma_1$
with $g_1 \neq e$. Since $\Gamma_1 \cap \Gamma_2 = \{ e \}$
we know that $g_1 \not\in \Gamma_2$ hence $g_1 g_2.x_1 \neq x_1$
and $c(x_2) = \alpha_1(g_1g_2.x_1)=c(g_1g_2.x_1)$. Moreover
$g_1 \neq e$ and $x_2 = g_2.x_1$ hence $g_1g_2.x_1 \neq x_2$. It
follows that $c(x_2) = c(g_1g_2.x_1) = \alpha_2(g_1g_2.x_1)
= \alpha_2(g_2g_1.x_1) = \alpha_2(g_1.x_1)$. We have
$g_1.x_1 \neq x_2 = g_2.x_1$ since $g_1 \not\in \Gamma_2$,
hence $c(x_2) = \alpha_2(g_1.x_1) = c(g_1.x_1) = \alpha_1(g_1.x_1) =
\alpha_1(x_1) < c(x_1)$, contradicting $c(x_1) < c(x_2)$.

As a fourth step we prove that there exists $(g_1,g_2)\in
\Gamma_1 \times \Gamma_2$ such that $g_1,g_2 \neq e$ and $x_2
\neq g_1 g_2.x_1$. We argue by contradiction. Let $g_2 \in \Gamma_2$
with $g_2 \neq e$. If, for all $g_1 \in \Gamma_1 \setminus \{ e \}$,
we have $g_1g_2.x_1 = x_2$ then $|\Gamma_1|=2$ by freeness of
the $\Gamma$-action. Similarly, we get $|\Gamma_2| = 2$. Since
$\Gamma$ is cyclic, $|\Gamma_1| = |\Gamma_2|$ implies $\Gamma_1 = \Gamma_2$
contradicting $\Gamma_1 \cap \Gamma_2 = \{ e \}$.

We can now conclude the proof. Let  $(g_1,g_2)\in
\Gamma_1 \times \Gamma_2$ such that $g_1,g_2 \neq e$ and $x_2
\neq g_1 g_2.x_1$. Then $\alpha_1(x_1) = \alpha_1(g_1.x_1) = c(g_1.x_1)$.
Moreover
$$
\begin{array}{lcll}
c(g_1.x_1) & = & \alpha_2(g_1.x_1) & \mbox{ since $x_2 \not\in \Gamma_1.x_1
\Leftrightarrow x_1 \not\in \Gamma_1.x_2$} \\
 & = & \alpha_2(g_2g_1.x_1) & = \alpha_2(g_1g_2.x_1) \\
 & = & c(g_1g_2.x_1) & \mbox{since $x_2 \neq g_1g_2.x_1$} \\
 & = & \alpha_1(g_1g_2.x_1) & \mbox{since $\Gamma_1 \cap \Gamma_2 = \{ e \}
\Rightarrow g_1 g_2.x_1 \neq x_1$} \\
 & = & \alpha_1(g_2.x_1) & = c(g_2.x_1) \mbox{ since $g_2\neq e$} \\
 & = & \alpha_2(g_2.x_1) & \mbox{ since $x_2\not\in \Gamma_2.x_1$}\\
 & = & \alpha_2(x_1) & = c(x_1) \mbox{ since $x_1\neq x_2$.} \\
\end{array}
$$
It follows that $\alpha_1(x_1) = c(x_1)$, contradicting $\alpha_1(x_1)
\neq c(x_1)$.
\end{proof}
\begin{figure}
\resizebox{5in}{!}{\includegraphics{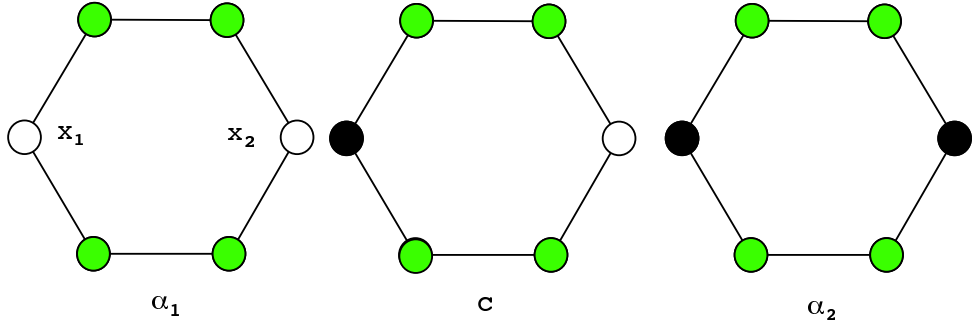}}
\caption{Necessity of assumption $\alpha < c$ instead of $\alpha \perp c$}
\label{neckfig6}
\end{figure}
\begin{figure}
\resizebox{5in}{!}{\includegraphics{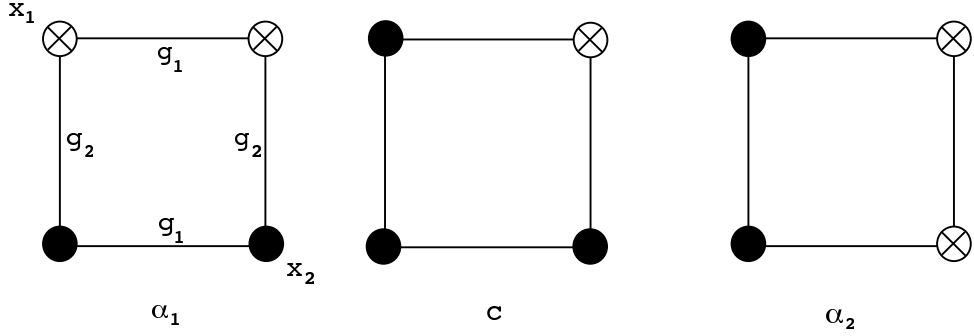}}
\caption{$\Gamma = < g_1,g_2 \ | \ g_1^2=g_2^2= e>$ and $\otimes < \bullet$}
\label{neckfig7}
\end{figure}

Figure \ref{neckfig6} illustrates the necessity of considering
necklaces with ordered pearls, and figure \ref{neckfig7}
shows that the assertion is false if $\Gamma$ is not cyclic.
However, the reader can check that the proof provided here works
for $\Gamma$ a (finite) commutative group with at most one
subgroup of order 2. 

\subsection{How many twins can one have ?}
Our goal is to study which necklaces appear
in pairs while fading one pearl in
a given necklace $c$.

\begin{lemma} \label{lemnecktwins}
Let $\alpha_1,\beta_1,\alpha_2,\beta_2 < c$ such that
$|\{\alpha_1,\beta_1,\alpha_2,\beta_2 \}| = 4$,
$\beta_1 = \gamma_1 .\alpha_1$, $\beta_2 = \gamma_2. \alpha_2$ with
$<\gamma_2 > \subset <\gamma_1 >$, $\gamma_1,\gamma_2 \not\in \Aut(c)$.
Then $\exists ! \mathcal{O} \in X/<\gamma_1>$ such that
$$
\left| c_{|\mathcal{O}}^{-1} \left( \max c(\mathcal{O}) \right) \right|
= |\mathcal{O}| - 1
$$
and, for all $P \in X/<\gamma_1>$,
$P \neq \mathcal{O} \Rightarrow |c(P)| = 1$.
\end{lemma}
\begin{proof}
Let $\Gamma_i = <\gamma_i>$. By lemma \ref{lemneckrestsym}, our assumptions imply
$\Aut(c) = 1$. Let $v_1,u_1,\mathcal{O}_1 = \mathcal{O}$
be given by lemma \ref{lemneckbasic} for $\gamma = \gamma_1$, and
$n = |\Gamma_1|$. In particular $\alpha_1(u_1) < c(u_1)$ and
$\beta_1(v_1) < c(v_1)$. For all $P \in X/\Gamma_1 \setminus \{ \mathcal{O} \}$
we have $|c(P)| = 1$ and $c(\mathcal{O}) = \{ \clubsuit, \spadesuit \}$
with $\clubsuit \neq \spadesuit$. We identify $< \gamma_1 >$ and
$\mathcal{O}$ with $\Z/(n)$,
$\gamma_1$ with $\bar{1}$, $v_1$ with $\bar{0}$ and
$u_1$ with $\bar{m}$ for some $m \in [0,n-1]$.
By assumption $u_1 = v_1$ hence $m \geq 1$, and we have $m \leq n-2$
because $\gamma_1 \not\in \Aut(c)$. We can assume $\clubsuit = c(v_1)$.
Among the statements of lemma \ref{lemneckbasic} we have
$|\beta_1(\OO)| = 2$ hence the set
$$
\beta_2(\OO) = \{ \beta_1(\bar{m}) , \beta_1(\bar{0}) , \beta_1(-\bar{1}) \}
= \{ c(\bar{m}), \beta_1(\bar{0}), c(-\bar{1}) \}
= \{ \clubsuit, \beta_1(v_1) , \spadesuit \}
$$
has cardinality 2. Since $\beta_1(\bar{0}) < c(\bar{0}) = \clubsuit$
it follows that $\beta_1(\bar{0}) = \spadesuit$ and $\spadesuit <
\clubsuit$, that is $\clubsuit = \max c(\OO)$. We then only
need to show that $m = n-2$.

Let $v_2,u_2,\mathcal{O}_2$
be given by lemma \ref{lemneckbasic} for $\gamma = \gamma_2$. Since
$\Gamma_2 \subset \Gamma_1$, the $\Gamma_2$-orbit
$\mathcal{O}_2$ is included in some $\Gamma_1$-orbit $P$.
But $|c(P) | \geq |c(\OO_2)| = 2$ and $c$ is constant
on every $\Gamma_1$-orbit besides $\OO$ hence $P = \OO$
and $\OO_2 \subset \OO$. In $\OO = \Z/(n)$ we identify
$v_2$ with $\bar{a}$ for some $a \in [0,n-1]$.

Assume first that $\Gamma_2 \varsubsetneq \Gamma_1$. Then
$\Gamma_2$ is generated by $\gamma_1^r$ for some
$r \geq 2$ dividing $n$. Since $|c(\OO_2)| = 2$
and $\OO_2 \subset \OO$ we have $c(\OO_2) = c(\OO)$
and $c(\bar{a}) = \max c(\OO) = \clubsuit$. It follows that
$\bar{a} \in [v_1,u_1]_{\gamma} = [\bar{0},\bar{m}]_{\bar{1}}$
and, since $m \geq 1$, there exists $b \in [0,n-1]$ such that
$\bar{b} \in \{ \overline{a-1},\overline{a+1} \}$
and $c(\bar{b}) = \clubsuit$. Since $\Gamma_1 \neq \Gamma_2$
and $\bar{a} \in \OO_2$ we have $\bar{b} \not\in \OO_2$.
Moreover $c$ is constant on every $\Gamma_2$-orbit different
from $\OO_2$ hence, for all $x \in \Z$,
the congruence $x \equiv b \mod r$  implies $c(\bar{x}) = \clubsuit$.
In particular there exists $x,x+r \in [0,n-1]$ such that
$c(\bar{x}) = c(\overline{x+r}) = \clubsuit$, hence
$c(\bar{z}) = \clubsuit$ for all $z \in [x,x+r]$. Now every
$\Gamma_2$-orbit in $\OO$ intersects $[\bar{x},\overline{x+r}]_{\bar{1}}$
hence $c(P) = \{ \clubsuit \}$ for all $P \in \OO/\Gamma_2$
with $P \neq \OO_2$ and $c(\OO \setminus \OO_2) = \{ \clubsuit\}$.
Let $\bar{x}, \bar{y} \in \mathcal{O}_2$
such that $c(\bar{x}) = x(\bar{y}) = \spadesuit$ with $x,y \in [0,n-1]$.
If $x \neq y$ we may assume $x < y$ hence
$m<x<y \leq n-1$ and $c(\overline{y-1}) = \spadesuit$
contradicting $\overline{y-1} \not\in \mathcal{O}_2$.
It follows that $x=y$, that is there is only one $x \in \mathcal{O}_2$ such that
$c(x) = \spadesuit$ and $|c^{-1}_{|\mathcal{O}}(\clubsuit)|
= |\mathcal{O}|-1$.

Now assume that $\Gamma_2 = \Gamma_1$, and let $r \in [0,n-1]$
with $\gcd(r,n) = 1$ such that $\gamma_2 = \gamma_1^r$.
Since $\Gamma_2 = \Gamma_1$ and $\OO_2 \subset \OO_1 = \OO$ we
have $\OO_2 = \OO$. In particular $c(v_2) = c(\bar{a})
= \max c(\OO) = \clubsuit$ and $|[v_2,u_2]_{\gamma_2}|
= |[v_1,u_1]_{\gamma_1}| =m+1$. Let $\varphi : x \mapsto
\bar{a} + r x \in \Aff_n$ (see section \ref{sectcyclic}
for the definition of $\Aff_n$). Since $\gamma_2 = \gamma_1^r$
we have $[v_2,u_2]_{\gamma_2} = \varphi([v_1,u_1]_{\gamma_1})=
\{ \varphi(\bar{0}), \dots ,\varphi(\bar{m}) \}$. On
the othe hand $[v_2,u_2]_{\gamma_2} = \{ x \in \OO \ | \ 
c(x) = \clubsuit \} = \{ \bar{0},\dots, \bar{m} \}$.
If $m < n-2$ lemma \ref{lemneckaffine} implies that
either $\varphi = \Id$ or $\forall x \ \varphi(x) = \bar{m}-x$.
The first case implies
$\alpha_1 = \alpha_2, \beta_1 = \beta_2$ and the second one
implies $\alpha_1 = \beta_2 , \alpha_2 = \beta_1$. Both thus
yield a contradiction, hence $m = n-2$ and 
$|c^{-1}_{|\mathcal{O}}(\clubsuit)|
= |\mathcal{O}|-1$.
\end{proof}

This lemma is a particular case of the following proposition,
and will be used in its proof.

\begin{prop} \label{propmulticompo2}
Let $\alpha_1,\beta_1,\alpha_2,\beta_2 < c$ such that
$|\{\alpha_1,\beta_1,\alpha_2,\beta_2 \}| = 4$,
$\beta_1 = \gamma_1 .\alpha_1$, $\beta_2 = \gamma_2. \alpha_2$ with
$\gamma_1,\gamma_2 \not\in \Aut(c)$. Let $\Gamma' = < \gamma_1,\gamma_2>$.
Then $\exists ! \mathcal{O} \in X/\Gamma'$ such that
$$
\left| c_{|\mathcal{O}}^{-1} \left( \max c(\mathcal{O}) \right) \right|
= |\mathcal{O}| - 1
$$
and, for all $P \in X/\Gamma'$,
$P \neq \mathcal{O} \Rightarrow |c(P)| = 1$. Moreover $c(x) = \alpha_i(x)
=\beta_i(x)$ for all $x \not\in \mathcal{O}$.
\end{prop}
\begin{proof}
Let $\Gamma_i = < \gamma_i>$. Lemma \ref{lemneckbasic} provides
$v_1,u_1,v_2,u_2,\OO_1,\OO_2$. We let $X' = \Gamma'.v_1$,
$\alpha'_i,\beta'_i, c'$ the restriction of $\alpha_i,
\beta_i,c$ to $X'$ and we denote by $\Aut(c')$
the stabilizer of $c'$ under the action of $\Gamma'$. Since
$\OO_1 \subset X'$, we know that $c$ is constant on
the $\Gamma_1$-orbits not included in $X'$, therefore
$\gamma_1 \not\in \Aut(c) \Rightarrow \gamma_1 \not\in \Aut(c')$.
Moreover $\OO_1 \subset X'$ implies $\alpha'_1 \perp c$,
$\beta'_1 \perp c$, $\beta'_1 = \gamma_1. \alpha'_1$,
whence $\Aut(c') = \{ e \}$ by lemma \ref{lemneckrestsym}. We show
that $\OO_2 \subset X'$. From $\beta'_2 = \gamma_2.\alpha'_2$ and
$\gamma_2 \neq 1$ we deduce $\alpha'_2 \neq c'$, because
otherwise we would have $\beta'_2 = \gamma_2.c'$ hence
\begin{itemize}
\item either $\beta'_2 = c'$ and $\gamma_2 \in \Aut(c') \setminus \{ 1 \}$,
contradicting $\Aut(c') = \{1\}$,
\item or $\exists ! x \in X' \ \beta'_2(x) = c'(x)$ contradicting
$\forall y \in Y \ |(\beta'_2)^{-1}(y)| = |c'^{-1}(y)|$ since
$\beta'_2 = \gamma_2.c'$.
\end{itemize}
Likewise, we have $\beta'_2 \neq c'$ hence $\alpha'_2 \perp c'$,
$\beta'_2 \perp c'$. Since $\beta'_2 = \gamma_2. \alpha'_2$
and $\gamma_2 \not\in \Aut(c') = \{ 1 \}$, lemma
\ref{lemneckbasic} implies that $X'$ contains some $\Gamma_2$-orbit
on which $c$ takes two distinct values, hence $\mathcal{O}_2 \subset X'$.

If $\Gamma_1 \subset \Gamma_2$ or $\Gamma_2 \subset \Gamma_1$, that
is $\Gamma' = \Gamma_1$ or $\Gamma' = \Gamma_2$,
lemma \ref{lemnecktwins} gives the conclusion so from now on we exclude
these cases. This assumption implies in particular
that there exists $P \in X'/\Gamma_1$, $Q \in X'/\Gamma_2$
with $P \neq \mathcal{O}_1$, $Q \neq \mathcal{O}_2$.
Since $\Gamma' = \Gamma_1 \Gamma_2$ acts freely and transitively
on $X'$, lemma \ref{lemneckclasses} implies that $c$ is constant
on both $X'\setminus \OO_1$ and $X'\setminus \OO_2$. Now
$\Gamma_1 \not\subset \Gamma_2$ and $\Gamma_2 \not\subset \Gamma_1$
implies that $\Gamma_1 \cup \Gamma_2$ is not
a subgroup of $\Gamma'$ hence $\Gamma_1 \cup \Gamma_2 \neq \Gamma'$.
By lemma \ref{lemneckclasses} there exists $v_0 \in \OO_1
\cap \OO_2$ hence $\OO_1 \cup \OO_2 = (\Gamma_1 \cup \Gamma_2).v_0$
and $|\OO_1 \cup \OO_2| = |\Gamma_1 \cup \Gamma_2 | < |\Gamma'|
= |X'$. It follows that $X' \neq \OO_1 \cup \OO_2$ therefore $c$ is
also constant on $X' \setminus (\OO_1 \cap \OO_2)$.
We let $\clubsuit \in Y$
denote this value taken by $c$. Since $Q \cap \mathcal{O}_1 \neq
\emptyset$ we have $\clubsuit \in c(\mathcal{O}_1)$ and
similarly $P \cap \mathcal{O}_2 \neq
\emptyset \Rightarrow \clubsuit \in c(\mathcal{O}_2)$. Let
$\spadesuit \in Y$ such that $c(\mathcal{O}_1) = \{ \clubsuit, \spadesuit \}$.
We have $\clubsuit \neq \spadesuit$ since $|c(\mathcal{O}_1)| = 2$.
Then $c(\mathcal{O}_1) = c(X') = c(\mathcal{O}_2)$ since $|c(\mathcal{O}_2)| = 2$.

If $c$ takes twice the value $\spadesuit$ on $\mathcal{O}_1$, by lemma
\ref{lemneckbasic} there exists $x \in \mathcal{O}_1$ such that
$c(x) = c(\gamma.x) = \spadesuit$, hence $x,\gamma_1.x \in \mathcal{O}_1 \cap \mathcal{O}_2$
and $\gamma_1 \in \Gamma_1 \cap \Gamma_2$ since
$\Gamma_1,\Gamma_2$ act freely transitively on $\mathcal{O}_1 \cap \mathcal{O}_2$.
In particular $\Gamma_1 \subset \Gamma_2$,
contradicting our assumption.

Since $c$ takes the value $\clubsuit$ on all the others $\Gamma_1$-orbits,
it follows that there exists a unique $x \in X'$ such that
$c(x) = \spadesuit$, and $]u_1,v_1[_{\gamma_1} = ]u_2,v_2[_{\gamma_2}
= \{ x \}$. In particular $v_1 = \gamma_1.x$ and $v_2 = \gamma_2.x$.
As in the proof of lemma \ref{lemnecktwins}, the existence
of $\alpha'_2 < c$ implies $\clubsuit > \spadesuit$. Letting
$\OO = X'$ we then have $|c_{|\OO}^{-1}(\max c(\OO))| = |\OO|-1$.
Let now $R \in X/\Gamma'$ with $R \neq \OO$.
Since $\mathcal{O}_1 \subset \mathcal{O}$
and $\mathcal{O}_2 \subset \mathcal{O}$ we know that $c$ is constant on
each $\Gamma_1$-orbit and each $\Gamma_2$-orbit in $R$. We deduce
from lemma \ref{lemneckclasses} that $|c(R)| = 1$ and the conclusion.
\end{proof}

\begin{figure}
\resizebox{5in}{!}{\includegraphics{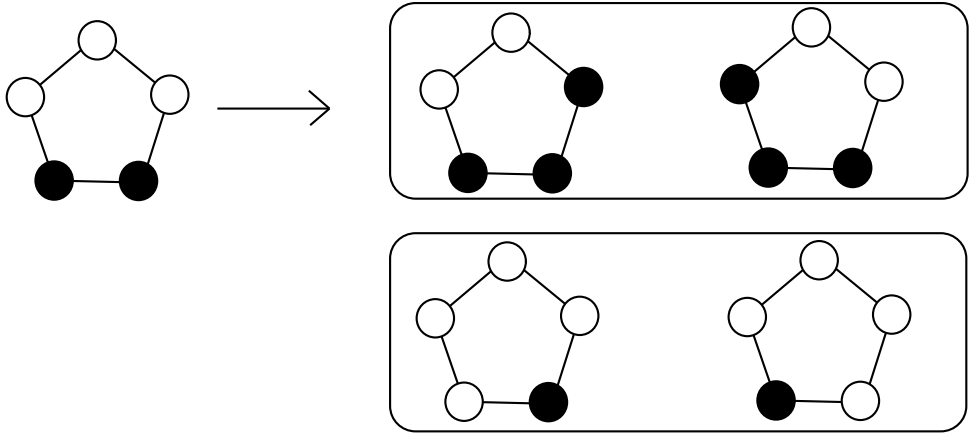}}
\caption{Necklaces with unordered pearls}
\label{neckfig3}
\end{figure}

\begin{figure}
\resizebox{5in}{!}{\includegraphics{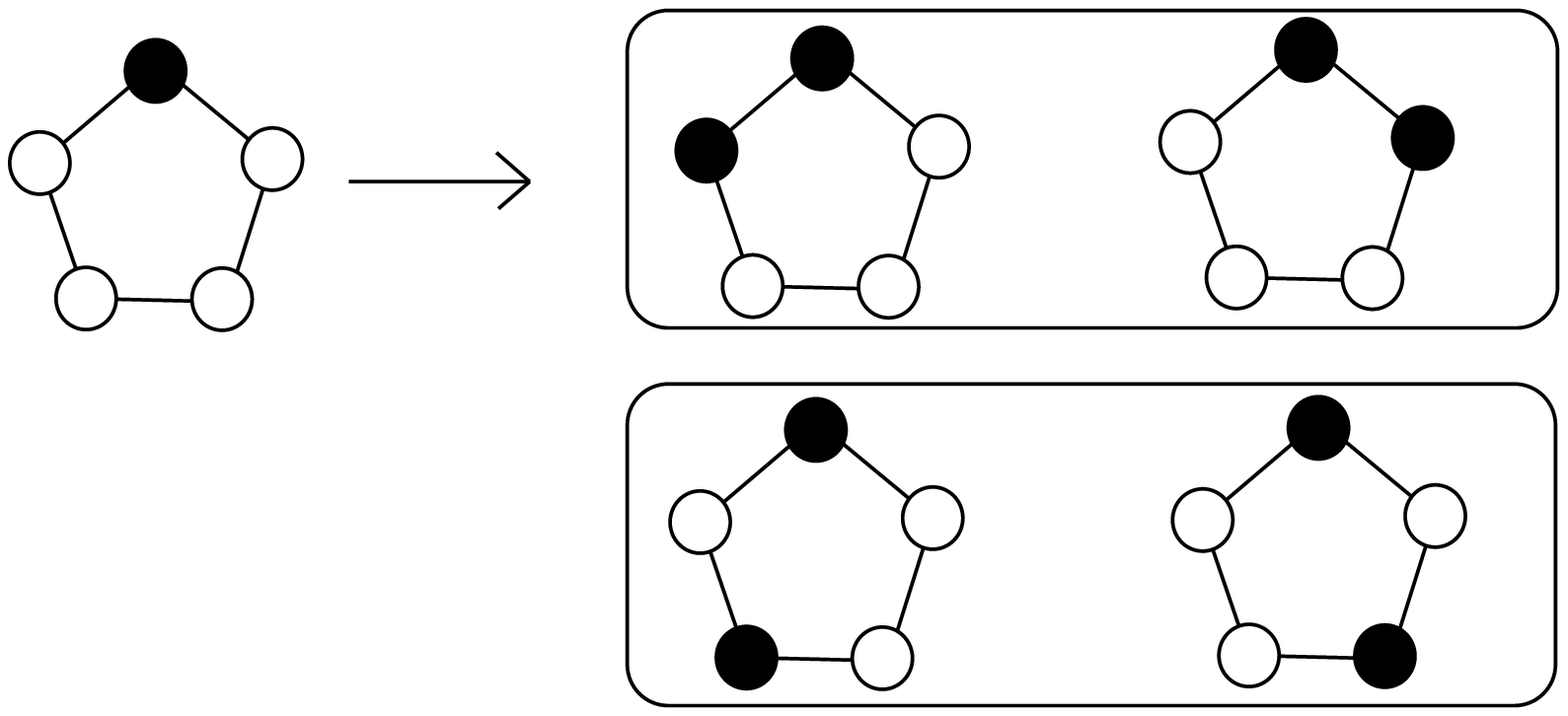}}
\caption{Necklace providing several twins (1)}
\label{neckfig4}
\end{figure}

\begin{figure}
\resizebox{5in}{!}{\includegraphics{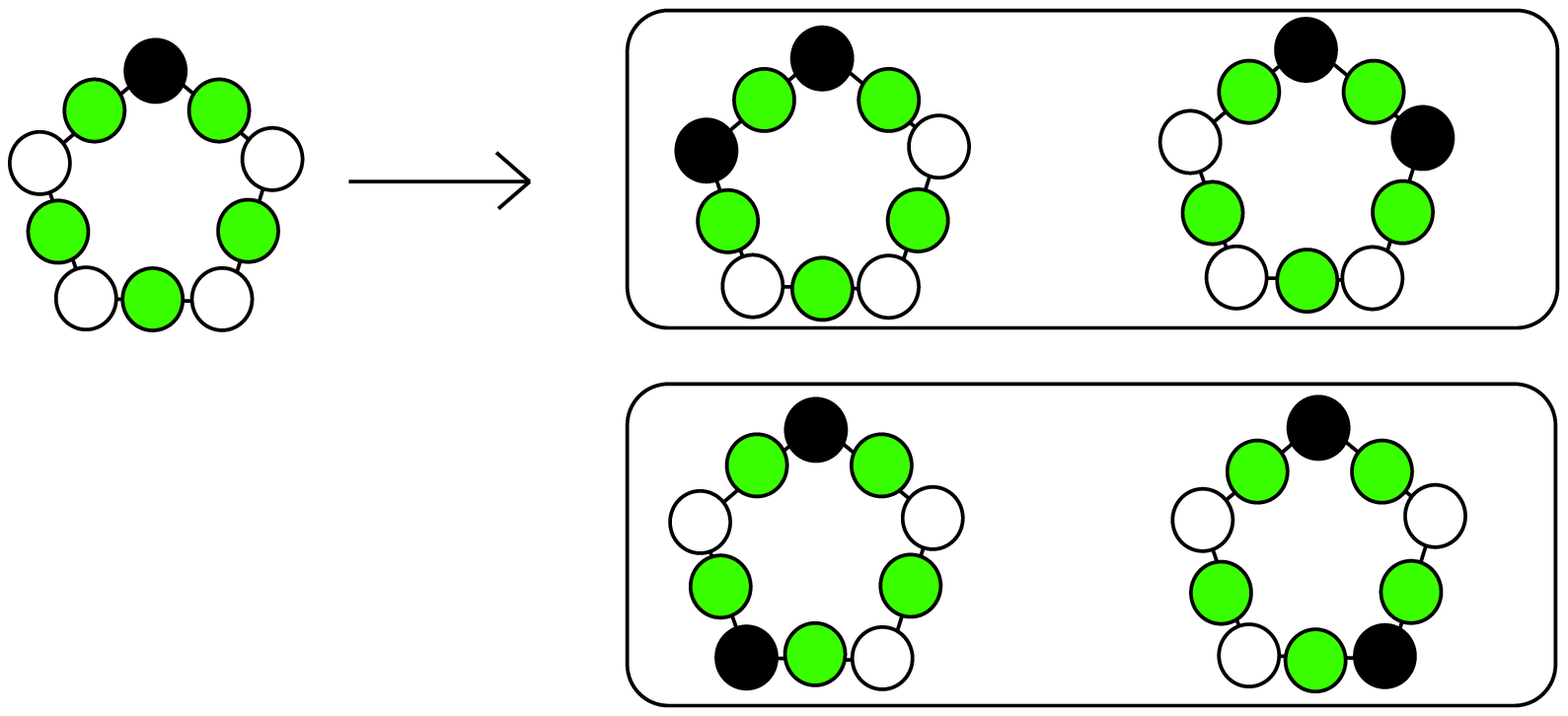}}
\caption{Necklace providing several twins (2)}
\label{neckfig5}
\end{figure}

Note that, for this proposition, we really need to
put a total order on $Y$. Figure \ref{neckfig3} shows a simple
necklace with four of its children which gives a counterexample
to the proposition if $<$ is replaced by the
weaker relation $\perp$. Figures \ref{neckfig4} and \ref{neckfig5} show
typical examples of necklaces (with ordered pearls, where
black is smaller than white) having several twins.

\begin{cor}
Under the assumptions of proposition \ref{propmulticompo2}, we have
$| \mathcal{O} | \geq 5$.
\end{cor}

\begin{proof}
Let $x \in \mathcal{O}$ such that $c(x) = \min c(\mathcal{O}) = \spadesuit$
and $\clubsuit = \max c(\mathcal{O})$. We have $\alpha_i(x) = c(x)$
because otherwise $c(x) \not\in \alpha_i(\mathcal{O})$. Since
$\beta_i = \gamma_i.\alpha_i$ and $\beta_i < c$ this implies
$\beta_i(x) = c(x)$ hence $c(x) \in \beta_i(\mathcal{O}) = 
\alpha_i(\mathcal{O})$, a contradiction. Similarly, $\beta_i(x) = c(x)$.

We also have, for $y \in \mathcal{O} \setminus \{ x \}$, $\beta_i(y)
\neq c(y) \Rightarrow \beta_i(y) = c(x)$. Indeed, $\beta_i(\gamma_i.x)
= \alpha_i(x) = c(x) \neq c(\gamma_i.x)$ by $\gamma_i \neq 1$
and the proposition. Since $\beta_i < c$ this implies $\gamma_i.x = y$
and $\beta_i(y) = c(x)$. Similarly, $\alpha_i(y) \neq c(y) \Rightarrow
\alpha_i(y) = c(x)$. In particular, $|\mathcal{O} \setminus \{ x \} |
\geq |\{ \alpha_1,\alpha_2,\beta_1,\beta_2 \} | = 4$ and the conclusion.
\end{proof}


\begin{table}
\caption{Table for exceptional groups}\label{tableexc}
$$
\begin{array}{|c|c|l|}
\hline
 \mbox{Group}&\mbox{Multiplicities}&\mbox{Generators of parabolic}\\
\hline
G_4&1&<s>\\
G_5&1&<s>\\
G_6&1&<t>\\
G_7& 1&
  <t>\\
G_8 & 1&<s>\\
G_9&1&
<t>\\
G_{10}&1&
<t>\\
G_{11}
  &1&
<u>\\
G_{12}&2&<s>\\
G_{13}&2&
<s>\\
G_{14}&2&
<s>\\
G_{15}&2&
<s>\\
G_{16}&2&
<s>\\
G_{17}&2&
<t>\\
G_{18}&2&
<s>\\
G_{19}&2& <t>\\
G_{20}&2&
<s>\\
G_{21}&2&
<t>\\
G_{22}&3&
<s>\\
G_{24}&2&<s,t>\\
G_{25}&1&
  <s,t>\\
G_{26}&1&
<s,t>\\
G_{27}&3&
  <t,u>\\
G_{29}&2&<s,t,u>\\
G_{31}&2&<s,t,u,v>\\
G_{32}&2&<s,t,u>\\
G_{33}&2&
<s,t,u,v>\\
G_{34}&2&<s,t,u,v,w>\\
\hline
\end{array}
$$
\end{table}

\vfill
\eject

\end{document}